\documentclass[preprint,3p,sort&compress]{elsarticle}

\journal{peer-review}
\usepackage[utf8]{inputenc}
\usepackage{bm}
\usepackage{amsmath}
\usepackage[hidelinks]{hyperref}
\hypersetup{colorlinks=true,linkcolor=blue,urlcolor=blue}
\usepackage{multicol}
\usepackage{booktabs}
\usepackage{amsthm}
\usepackage[table,xcdraw,dvipsnames]{xcolor}
\usepackage{graphicx}
\usepackage{float}
\usepackage{subcaption}
\usepackage{pdflscape}
\usepackage{array}
\usepackage{multirow}
\usepackage{cleveref}
\usepackage{amsfonts}
\usepackage{amssymb}
\usepackage{fourier-orns}
\usepackage[labelfont=bf,labelsep=period]{caption}
\usepackage[labelfont=rm,font=footnotesize]{subcaption}

\numberwithin{equation}{section}

\DeclareMathAlphabet{\altmathcal}{OMS}{cmsy}{m}{n}
\DeclareMathAlphabet{\altmathcalb}{OMS}{cmsy}{b}{n}
\DeclareMathAlphabet{\mathcalboondox}{U}{BOONDOX-calo}{m}{n}
\DeclareMathAlphabet{\mathbbmsl}{U}{bbm}{m}{sl}
\newcommand{\Ir}{\altmathcal{I}}
\newcommand{\Sr}{\altmathcal{S}}
\newcommand{\Or}{\altmathcal{O}}
\newcommand{\Ur}{\altmathcal{U}}
\newcommand{\Vr}{\altmathcal{V}}
\newcommand{\Ar}{\altmathcal{A}}
\newcommand{\Br}{\altmathcal{B}}
\newcommand{\ur}{\mathcalboondox{u}}
\newcommand{\vr}{\mathcalboondox{v}}
\renewcommand{\wr}{\mathcalboondox{w}}
\newcommand{\ui}{\mathbbmsl{u}}
\newcommand{\vi}{\mathbbmsl{v}}
\newcommand{\Vi}{\mathbbmsl{V}}

\newcommand{\R}{\mathbb{R}}             
\newcommand{\vphi}{\varphi}				
\newcommand{\eps}{\varepsilon}	
\newcommand{\sg}{\sigma}	
\newcommand{\Gm}{\Gamma}
\newcommand{\Om}{\Omega}				
\newcommand{\lm}{\lambda}				
\newcommand{\N}{N} 				 		
\renewcommand{\ne}{{n_e}} 				
\newcommand{\B}[2]{v_{#1}^{\.#2}}		

\renewcommand{\.}{\!\:}									
\newcommand{\norm}[1]{\left\lVert#1\right\rVert} 		

\newdefinition{remark}{Remark}

\newcommand{\orcid}[1]{\href{https://orcid.org/#1}{\texorpdfstring{\includegraphics*[width=8pt]{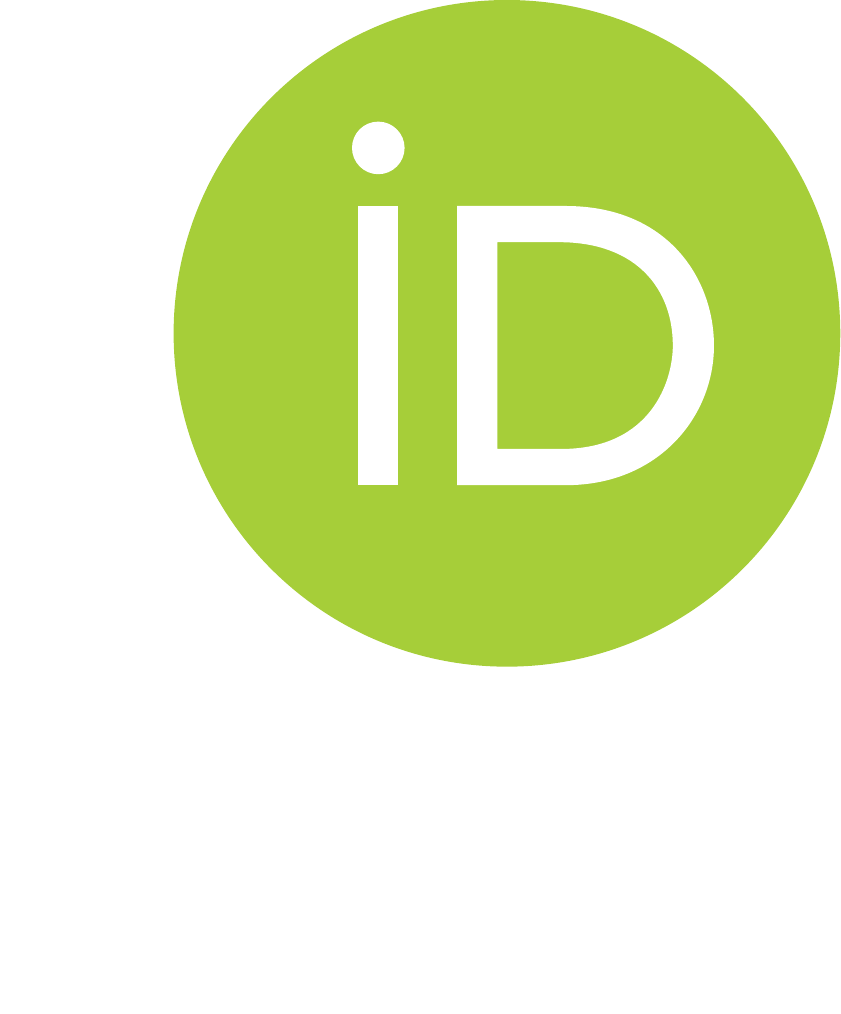}}~~}{}}

\definecolor{drot}{rgb}{0.7,0,0.1}

\newcommand{\az}[1]{{\color{blue}#1}}

\begin{document}
\baselineskip14pt
\sloppy

\begin{frontmatter}

\title{Machine Learning Discovery of Optimal Quadrature Rules for \\Isogeometric Analysis}

\author[bcam]{Tomas Teijeiro\orcid{0000-0002-2175-7382}\corref{cor1}}
\cortext[cor1]{Corresponding author}
\ead{tteijeiro@bcamath.org}
\author[cunef]{Jamie M. Taylor\orcid{0000-0002-5423-828X}}
\author[bcam]{Ali Hashemian\orcid{0000-0002-8190-222X}}
\author[upv,bcam,ikerbasque]{David Pardo\orcid{0000-0002-1101-2248}}

\address[bcam]{{BCAM -- Basque Center for Applied Mathematics}, {Bilbao, Basque Country, Spain}}
\address[cunef]{{Department of Quantitative Methods, CUNEF University}, {Madrid, Spain}}
\address[upv]{{University of the Basque Country (UPV/EHU)}, {Leioa, Basque Country, Spain}}
\address[ikerbasque]{{Ikerbasque -- Basque Foundation for Sciences}, {Bilbao, Basque Country, Spain}}

\begin{abstract}
We propose the use of machine learning techniques to find optimal quadrature rules for the construction of stiffness and mass matrices in isogeometric analysis~(IGA).
We initially consider 1D spline spaces of arbitrary degree spanned over uniform and non-uniform knot sequences, and then the generated optimal rules are used for integration over higher-dimensional spaces using tensor product sense.
The quadrature rule search is posed as an optimization problem and solved by a machine learning strategy based on gradient-descent. However, since the optimization space is highly non-convex, the success of the search strongly depends on the number of quadrature points and the parameter initialization.
Thus, we use a dynamic programming strategy that initializes the parameters from the optimal solution over the spline space with a lower number of knots. 
With this method, we found optimal quadrature rules for spline spaces when using IGA discretizations with up to~50 uniform elements and polynomial degrees up to 8, showing the generality of the approach in this scenario.
For non-uniform partitions, the method also finds an optimal rule in a reasonable number of test cases. 
We also assess the generated optimal rules in two practical case studies, namely, the eigenvalue problem of the Laplace operator and the eigenfrequency analysis of freeform curved beams, where the latter problem shows the applicability of the method to curved geometries.
In particular, the proposed method results in savings with respect to traditional Gaussian integration of up to 44\% in 1D, 68\% in 2D, and 82\% in 3D spaces.
\end{abstract}


\end{frontmatter}



\section{Introduction}
\label{sec:Introduction}

Developing efficient and accurate integration methods plays a crucial role in many numerical analysis techniques.
In the context of the isogeometric analysis (IGA)~\cite{Hughes2005}, a common practical approach for numerical integration is to use an element-wise Gaussian (EWG) quadrature rule when constructing system matrices in the sense of Galerkin discretizations. This follows the classical system construction technique of the finite element analysis (FEA).
However, it is known that there exist optimal quadrature rules for spline spaces of higher continuities, thus requiring a significantly fewer number of quadrature points than the classical EWG.
For instance, one may consider the fast matrix formation technique by~\citet{Calabro2017}, where each row of the system matrices is integrated by its own quadrature rule obtained by a linear system of equations.
\citet{Barton2020} propose weighted Gaussian quadrature rules for B-splines that require the minimum number of quadrature points while guaranteeing the exactness of integration with respect to the weight function.
Other techniques by~\citet{barton_gaussian_2016,bartonOptimalQuadratureRules2016, Barton2017} use a polynomial homotopy continuation~(PHC), a numerical scheme for solving polynomial systems of equations~\cite{Sommese2005}, to generate Gaussian quadrature rules for spline spaces of higher continuities.
To generate a Gaussian rule in a \textit{target} spline space, they built an associated \textit{source} space with known quadratures (e.g., a union of polynomial Gaussian rules) and transform the rule from the source space to the target space, while preserving the optimality.
These rules, in the limit (when the number of elements are very large), converge to the half-point rules of~\citet{hughes_efficient_2010}.
More research works on obtaining optimal or nearly optimal quadrature rules for IGA discretizations can be found in, e.g.,~\cite{ auricchio_simple_2012, Schillinger2014, Hiemstra2017, Barendrecht2018, Zou2021, Giannelli2022}.

On the other hand, the use of machine learning techniques to address numerical integration problems is gaining momentum in recent years. While the practical state-of-the-art for calculating approximate integrals with bounded errors are still Monte Carlo methods~\cite{MonteCarloIntegration2004}, kernel-based quadrature approaches have proved interesting convergence properties, making them suitable for high-dimensional problems, or when the integrand function is expensive to evaluate~\cite{NIPS2016_81c650ca, NEURIPS2018_6e923226}. Another example are Bayesian quadrature methods~\mbox{\cite{ ohaganBayesHermiteQuadrature1991, karvonenClassicalQuadratureRules2017, NEURIPS2019_165a59f7}}, which aim at providing a statistically accurate integration while minimizing the required number of evaluations of the integrand. Quadrature techniques are also a relevant stage on deep learning-based  approximations, for which no proper rules exist (see, e.g.,~\cite{riveraQuadratureRulesSolving2022}). All of the mentioned approaches try to calculate an approximate value of the integral, while the main contribution of this work lies in providing optimal quadrature rules, in the sense of calculating the exact integral with the minimum possible number of quadrature points.

In order to guarantee the optimality of the discovered rules, we first define the target number of quadrature points according to the half-point rule~\cite{hughes_efficient_2010}, and then set as trainable parameters the location and weight of each point. The loss function that completes the definition of the machine learning problem is based on the integration error of each individual basis function in the fixed finite dimensional spline space. However, this problem is highly non-convex, and the convergence of the learning process to the global optimum depends fundamentally on the initial value of the parameters~\cite{narkhedeReviewWeightInitialization2022}. For this reason, and exploiting the regularity of spline spaces with the same degree and continuity, but different numbers of elements, we devise a dynamic programming strategy~\cite{bellmanDynamicProgramming1966} that initializes the parameters of a given problem from the optimal solution of a simpler problem. The experimental results show that this strategy achieves general convergence on partitions with uniform elements.

The structure of the remainder of this paper is as follows:
\autoref{sec:ProblemFormulation} formally defines the target function spaces and the integration problem.
In \autoref{sec:Method}, we describe the developed optimization method, including the loss function and a dynamic programming strategy for parameter initialization.
Then, \autoref{sec:Results} shows the numerical results of the obtained quadrature rules for general spline spaces, their comparison with classical EWG rules, and the computational requirements. 
In \autoref{sec:CaseStudies}, we discuss the IGA results of two applied case studies.  
Finally, \autoref{sec:Conclusions} draws some conclusions and indicates possible directions for future research.


\section{Problem formulation}
\label{sec:ProblemFormulation}

\subsection{IGA discretization}
\label{sub:IGA}

\begin{figure}[h!]
     \centering
         \includegraphics[]{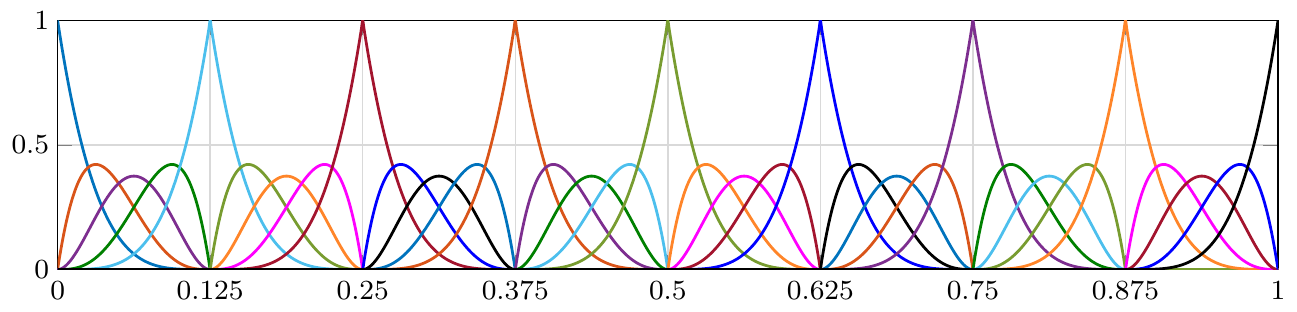}
        \caption{Basis functions of an example B-spline space with $p=4$, $c=0$, and $\ne=8$.}
        \label{fig:basis_examples}
\end{figure}

For the sake of simplicity, we consider an IGA discretization over the one-dimensional parameter space. The extension to the higher-dimensional spaces is trivial by taking tensor products.
Without loss of generality, we introduce the computational domain $\Om$ as a partition of $\ne$ elements in ${[0,1]}$, described by a sequence of non-repeating knots ${U:=\{u_j\}_{j=0}^{\ne}}$ such that ${0=u_0<u_1<\ldots<u_{\ne-1}<u_{\ne}=1}$.
We define the B-spline space $\Sr^c_p(U)$ of piecewise polynomial functions of degree $p$  with continuity $c$ at every {\em interior} knot in $U$, where ${p>c\geq 0}$. 
Let us consider the \textit{clamped} knot vector ${\Xi\supset U}$ as 
\begin{align}
\Xi:=\{\.\underbrace{0,\ldots,0}_{p+1},\,\underbrace{u_1,\ldots,u_1}_\mu,\,\underbrace{u_2,\ldots,u_2}_\mu\,\ldots,\,\underbrace{1,\ldots,1}_{p+1}\}=\{\xi_0, \xi_1, \ldots , \xi_{n+p+1}   \}\.,
\end{align}
where ${\mu:=p-c}$ is the multiplicity of interior knots and ${n+1}$ is the dimension of $\Sr^c_p(U)$.
We obtain the basis functions ${\{v_i^{\.p}\}_{i=0}^n}$ of this B-spline space by the Cox--De~Boor recursion~\cite{TheNURBSBook}. For a given parameter ${x\in\Om}$, it reads
\begin{align}
	\B{i}{0}(x)  &= 
	\begin{cases}\begin{aligned}
	&1\,,&&\xi_i\leq x <\xi_{i+1}\,, \\ \label{eq.BasisFunBi0} 
	&0\,,&&{\rm otherwise}\,,
	\end{aligned}\end{cases} \\
	\B{i}{p}(x) &= \dfrac{x-\xi_i}{\xi_{i+p}-\xi_i}\B{i}{p-1}(x)+
	\dfrac{\xi_{i+p+1}-x}{\xi_{i+p+1}-\xi_{i+1}}\B{i+1}{p-1}(x)\,. \label{eq.BasisFunBip}
\end{align}
When evaluating basis functions at a desired $x$, we find the corresponding nonzero knot span in~\eqref{eq.BasisFunBi0} and efficiently evaluate~\eqref{eq.BasisFunBip} avoiding any division by zero in dealing with repetitive knots. More details are given in \mbox{\cite[Algorithms A2.1 and A2.2]{TheNURBSBook}} and~\cite[Algorithms A1 and A2]{Hashemian2022}.
%
%
%
%
%
\autoref{fig:basis_examples} shows an example of B-spline bases characterized by $p=4$ and $c=0$ spanned over a uniform partition with $\ne=8$ elements.

\begin{remark}
We note that $\Sr_c^p(U)$ admits different basis functions. For example, an alternative set of bases are $\{x^i\}^p_{i=0}$ and ${\bigcup_{i=0}^{n_e-1}\{\text{ReLU}(x-u_i)^r\}^p_{r=c+1}}$. However, they produce excessively large condition numbers that lead to numerical instabilities when searching for the optimal quadrature rules.
\end{remark}


\subsection{Numerical integration}
\label{sub:NumericalIntegration}

Without loss of generality, we aim at finding optimal quadrature rules to exactly compute the following integrals when constructing the stiffness $K$ and mass $M$ matrices arising from the discretization of a partial differential equation (PDE):
\begin{alignat}{3}
K_{ij}&=\int_{\Om} v'_i \. v'_j \,d\Om \,, 	&\qquad i&=0,1,\ldots, n \,, \label{eq:Stiffness} \\
M_{ij}&=\int_{\Om} v_i \. v_j \,d\Om \,, 	&\qquad i&=0,1,\ldots, n \,, \label{eq:Mass}
\end{alignat}
where the prime denotes the derivative of the basis functions. 
For the sake of brevity, herein and in the following, we omit the superscript $p$ in referring to basis functions. We note that the integrals~\eqref{eq:Stiffness} and~\eqref{eq:Mass} are formed without taking into account the coefficients of the investigated PDE and the problem geometry. This is a common assumption when seeking to compute the optimal quadrature rules (see, e.g.,~\cite[Section~3.2]{hughes_efficient_2010} and~\cite[Section~2]{auricchio_simple_2012}).

\begin{remark}
Let $\Sr_c^p(U)$ be a spline space used for IGA discretization, the $L^2$ product of basis function derivatives in~\eqref{eq:Stiffness} entails integrating the space $\Sr_{c-1}^{2p-2}(U)$, while in \eqref{eq:Mass}, we integrate the space $\Sr_{c}^{2p}(U)$.
As a result, the minimal space that contains both stiffness and mass terms is $\Sr_{c-1}^{2p}(U)$. Furthermore, when using typical maximum-continuity IGA discretizations, $c$ is equal to ${p-1}$.
Thus, for brevity, we consider the space $\Sr_k^d(U)$ for the integrand, where $d$ is equal to $2p$\., and $k$ is up to $p-2$.
\label{rem:Spaces}
\end{remark}

We wish to find quadrature points ${X := \{x_i\}_{i=1}^q}$ and weights ${W := \{w_i\}_{i=1}^q}$ presuming that a given integral bilinear form ${B:\Sr_k^d(U)\to\R}$ is evaluated exactly on $\Sr_k^d(U)$. Thus, for any $f,g\in \Sr_c^p(U)$, we have
\begin{equation}
B(f,g) :=\sum\limits_{i=1}^q f(x_i)g(x_i)w_i\,.
\end{equation}
We note, however, that it suffices to produce an appropriate quadrature rule in a larger space for integrals, as we have that if ${f,g\in \Sr_c^p(U)}$, then ${f\cdot g\in \Sr_k^{d}(U)}$ in the sense of \autoref{rem:Spaces}. Thus, we aim to find quadrature rules for integration of the form
\begin{align}
I(f;X,W):=\sum\limits_{i=1}^q f(x_i)w_i\,,
\end{align}
for all $f\in \Sr_k^{d}(U)$, with the following constraints:
\begin{equation}
\begin{aligned}
   x_i, w_i &\in [0, 1]\,,\quad i=1,2,\ldots,q\,,\\
   \sum\limits_{i=1}^q w_i &= 1\,.
\label{eq:constraints}
\end{aligned}
\end{equation}

When constructing system matrices in IGA, a common practical approach is to use the classical element-wise Gaussian (EWG) quadrature rule that entails a ${(p+1)}$-point rule within each element ${\Om^e_j:[u_{j-1},u_j]}$, ${j=1,\ldots,\ne}$, for integrals~\eqref{eq:Stiffness} and~\eqref{eq:Mass}.
This element-wise rule is optimal only when using a discontinuous space (i.e., with continuity ${k=-1}$ at all elements interfaces), which could be the case in the classical FEA discretizations.
Thus, the extra cost of EWG rules in IGA may lead to significant computational overheads, particularly with higher-order discretizations and in higher-dimensional problems. 
For spline spaces, there always exists an optimal quadrature rule where the number of quadrature points $q$ is given by:
\begin{align}
d+l+1=2q\,,
\label{eq:OtimalQuadPoints}
\end{align}
where $l$ is the total number of interior knots, including their multiplicities, i.e., ${l := (\ne-1)\cdot(d-k)}$ assuming a uniform continuity $k$ at all interior knots. If ${d\ll l}$, then ${q\approx l/2}$, denoting that in the limit, when ${\ne\rightarrow\infty}$, the mentioned spline rules converge to the half-point rules of~\citet{ hughes_efficient_2010}, which are exact over the domain.
Therefore, we consider a quadrature rule to be optimal if $q$ is the minimum required number of quadrature points, while guaranteeing the exactness of the integration.

\begin{remark}
Based on some combinations of $p$, $c$, and $\ne$, the left-hand side in~\eqref{eq:OtimalQuadPoints} may become odd; thus, we obtain the optimal number of quadrature points by ${q=\lceil\frac{d+l+1}{2}\rceil}$. In such occasions, one may consider a higher space $\Sr_k^{d+1}(U)$ to have a unique quadrature rule following the idea of integrating over odd-degree spline spaces described in, e.g.,~\cite{bartonOptimalQuadratureRules2016}. However, this implies the over-integration of system matrices.
\label{rem:OddDegree}
\end{remark}


\section{Optimization method}
\label{sec:Method}

\subsection{Parameter space and loss function}
\label{sub:loss_function}
As $\Sr_k^{d}(U)$ is a finite-dimensional space, all norms upon it are equivalent, and we consider it to be equipped with the $L^2$ norm. Then, for each fixed $X,W$, we can view $I(\cdot;X,W)$ as a linear functional on $\Sr_k^{d}(U)$. In this case, the natural way to quantify the integration error would be using the dual norm defined over $\Sr_k^{d}(U)$. Whilst dual norms are generally difficult to evaluate, as we consider $L^2$ to be the norm on $\Sr_k^{d}(U)$, we have a straightforward numerically tractable expression for the induced dual norm. Specifically, if ${\Ir(f):=\int_\Om f(x)\,dx}$ is the exact integral of $f$, we have

\begin{equation}
\label{eq:dual_norm}
||\Ir-I(\cdot;X,W)||:=\max\limits_{f\in \Sr_k^{d}(U)}\frac{\left|\Ir(f)-I(f;X,W)\right|}{||f||_{L^2}}\,.
\end{equation}

Given any basis $\{v_i\}_{i=0}^n$ of $\Sr_k^{d}(U)$, we consider the Gram (a.k.a mass) matrix given by~\eqref{eq:Mass} and, by utilizing the Hilbert space structure of $L^2$ and its dual space, we define our loss function as:
\begin{equation}
\label{eq:loss_function}
L(X, W) := \max\limits_{f\in \Sr_k^{d}(U)}\frac{\left|e(f)\right|}{||f||_{L^2}}=\left(\sum\limits_{i,j=0}^n (M^{-1})_{ij}\.e(v_i)\.e(v_j)\right)^\frac{1}{2},
\end{equation}
where $e(v_i)$ is the integration error for the $v_i$ basis function, that is for a fixed $X$ and $W$:
\begin{equation}
e(v_i) := \left|\.\Ir(v_i) - I(v_i;X,W)\right|.
\end{equation}
If the basis were orthonormal, then $M$ and its inverse would be the identity matrix, and generally $M^{-1}$ is utilized to ensure consistency between the loss and the $L^2$ norm within a general basis.

\begin{figure}[!h]
     \centering
     \begin{subfigure}[t]{0.34\textwidth}
         \centering
         \includegraphics[width=\textwidth]{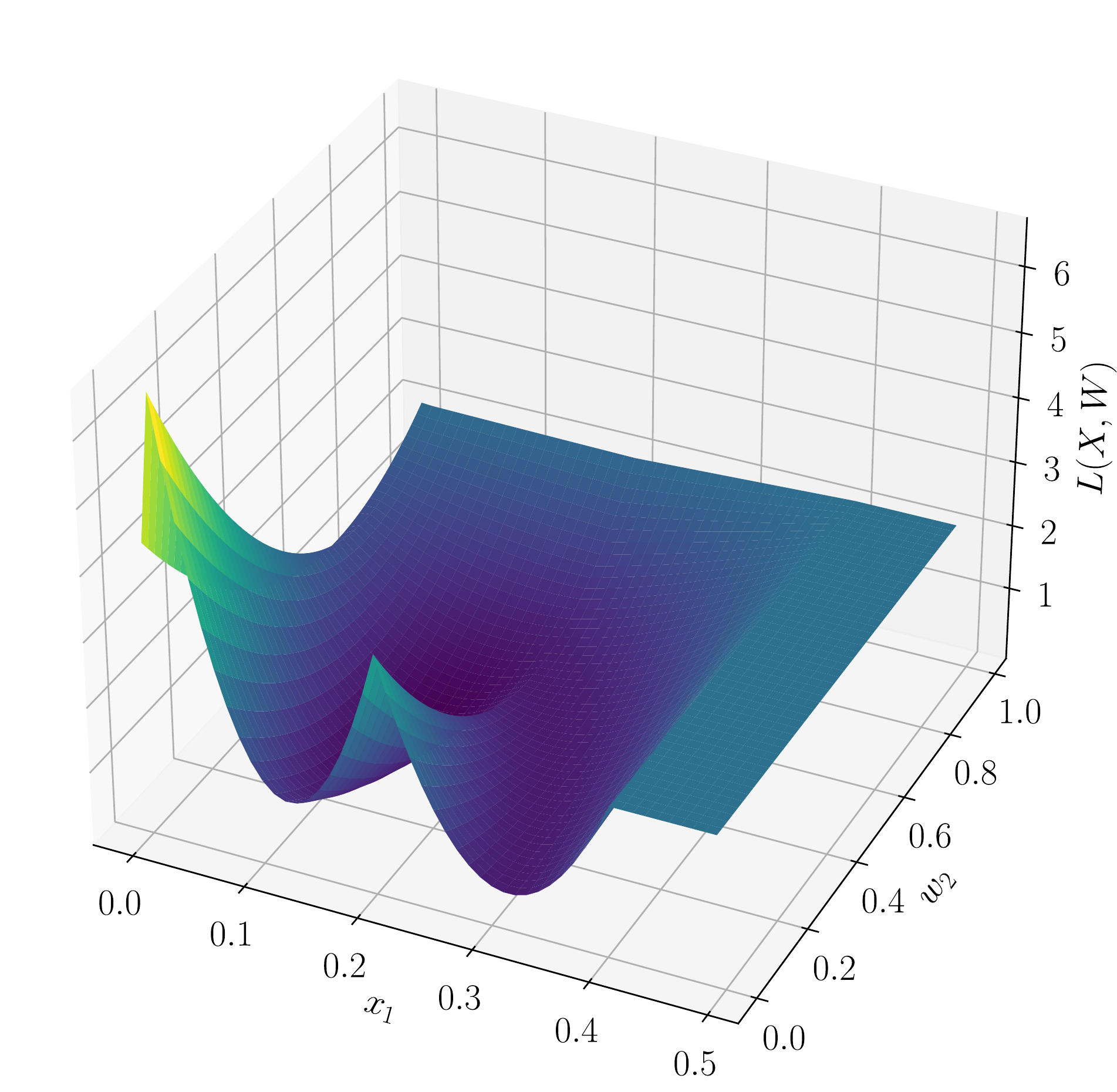}
         \caption{Original search space.}
         \label{fig:loss}
     \end{subfigure}
     \hfill
     \begin{subfigure}[t]{0.34\textwidth}
         \centering
         \includegraphics[width=\textwidth]{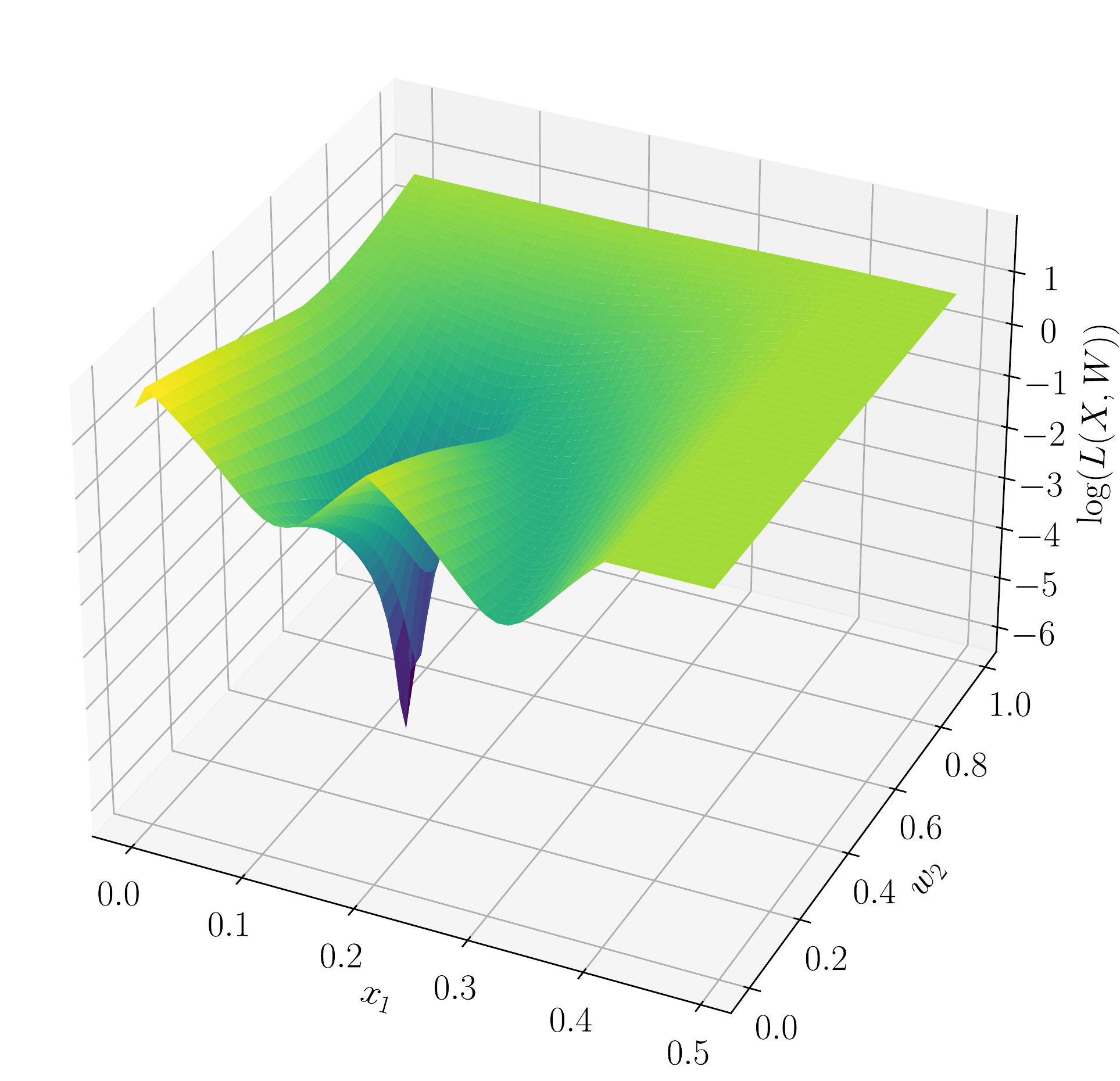}
         \caption{Logarithmic scale to highlight the single optimum.}
         \label{fig:logloss}
     \end{subfigure}
     \hfill
     \begin{subfigure}[t]{0.30\textwidth}
         \centering
         \includegraphics[width=\textwidth]{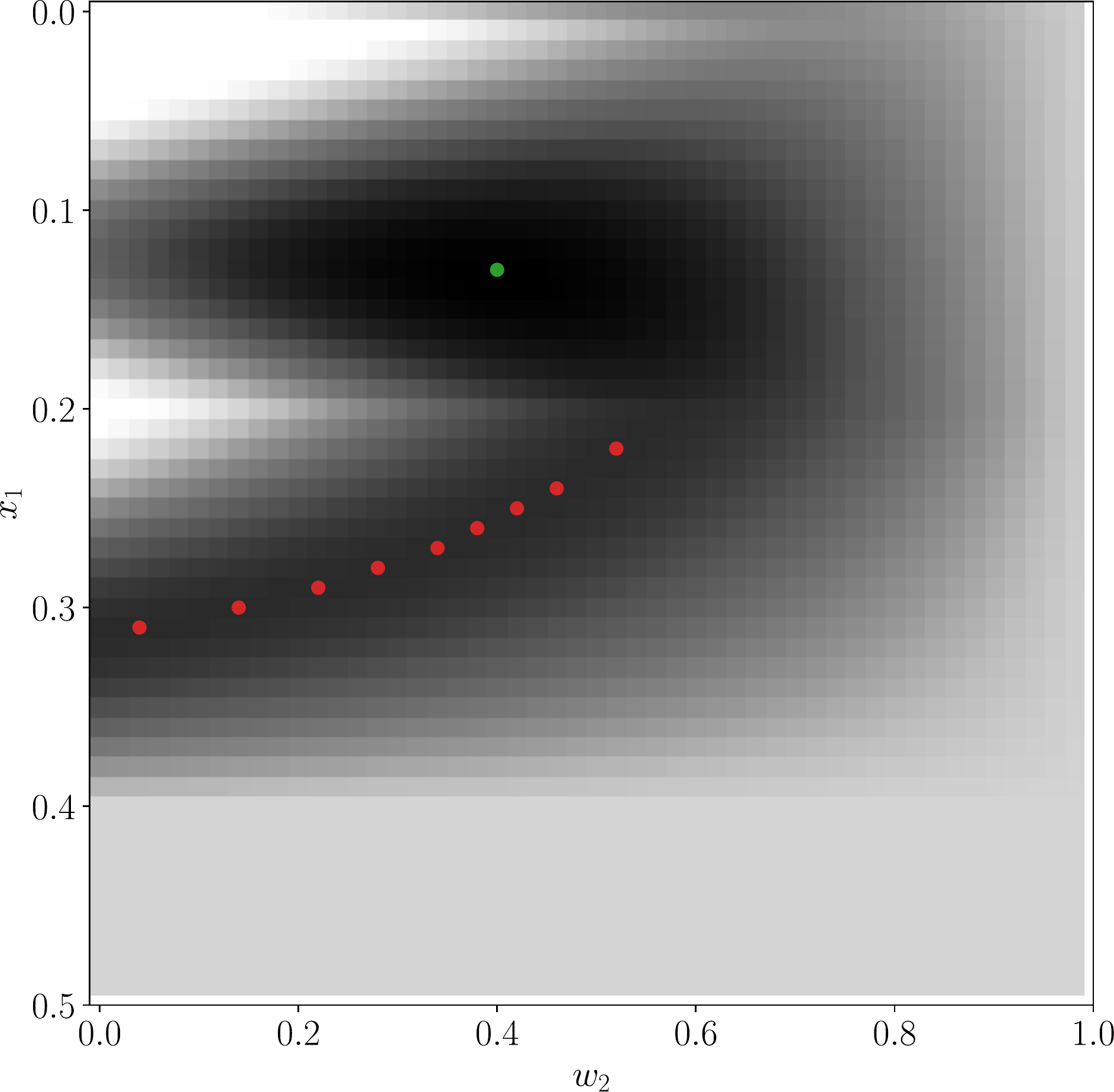}
         \caption{2D bird's eye view of the search space. The global optimum is marked in green, and local minima in red.}
         \label{fig:loss_2d}
     \end{subfigure}
        \caption{Illustrative example of the optimization space for $d=1$, $k=0$ and $\ne=5$.}
        \label{fig:search_space}
\end{figure}

In general, for a fixed number of quadrature points $q$, there are $2q$ unknown parameters corresponding to the point locations $\{x_i\}_{i=1}^q$ and weights $\{w_i\}_{i=1}^q$, and the optimization problem consists of minimizing $L(X,W)$ while respecting the constraints defined in~\eqref{eq:constraints}. Since \eqref{eq:loss_function} is differentiable, we approach this problem first by setting the optimal number of points to calculate as ${q=\lceil \frac{n}{2}\rceil}$, c.f.~\cite{hughes_efficient_2010}, and then by using a standard machine learning strategy based on gradient-descent. However, we observe that the optimization space is highly non-convex, and thus the search usually stagnates at local minima, being highly dependent on $q$ and the initial values of $x_i$ and $w_i$. This is illustrated in \autoref{fig:search_space}, where we show the optimization space for a simple problem ($d=1$, $k=0$ and ${\ne=5}$). In this space, $q$ is equal to three; but due to the symmetry constraints explained below~\eqref{eq:symmetry_constraints}, together with the basic constraints in~\eqref{eq:constraints}, the loss is determined by $x_1$ and $w_2$ only. We can see that even in this simple space there are multiple local minima, and the ability of gradient descent to effectively find the solution depends on the starting point for searching.

\subsection{Dynamic programming for parameter initialization}

Given the relevance of the initial values of the searching parameters for the convergence of the machine learning optimization, we devise a method based on dynamic programming~\cite{ bellmanDynamicProgramming1966} to initialize the parameter values in a fully deterministic way. In short, the initial parameter values for a given problem will be calculated from the optimal quadrature rule of a \textit{simpler} problem. This allows us to break the solution of a problem into a chain of simpler problems and to solve it recursively. Let us start by discussing how to solve problems with uniform elements, and then, we will build on those results to address more general problems with non-uniform elements.

\subsubsection{Uniform elements}
\label{sub:initialization_uniform}

\begin{figure}[!h]
     \centering
     \begin{subfigure}[b]{0.48\textwidth}
         \centering
         \includegraphics[width=\textwidth]{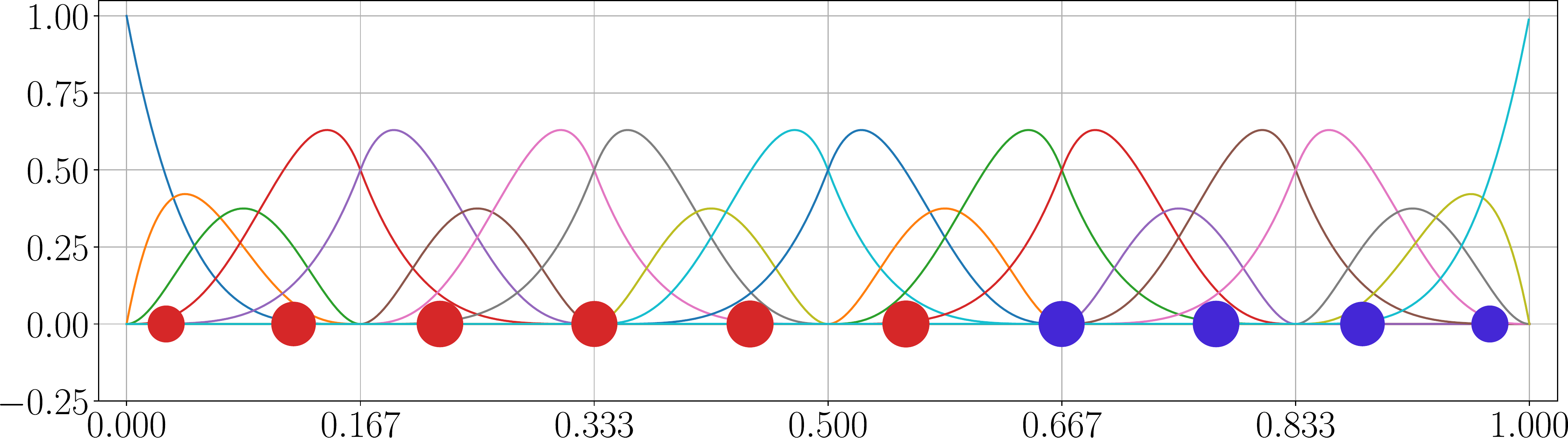}
         \caption{Optimal quadrature points for $d=4, k=1, \ne=6$.}
         \label{fig:4_1_6}
     \end{subfigure}
     \hfill
     \begin{subfigure}[b]{0.48\textwidth}
         \centering
         \includegraphics[width=\textwidth]{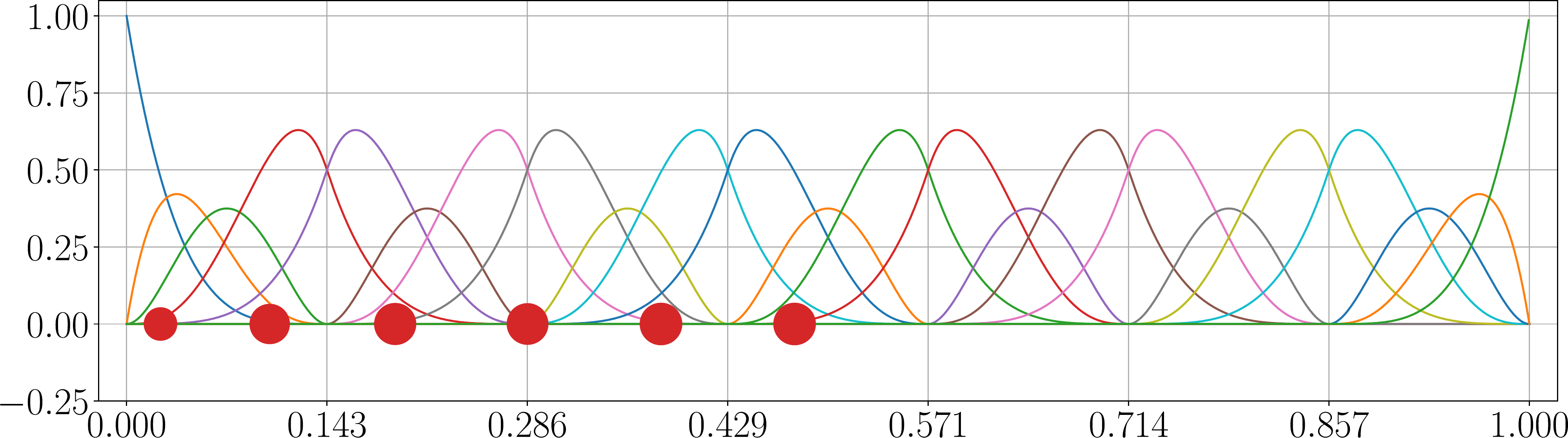}
         \caption{Initial parameter values for $d=4, k=1, \ne=7$.}
         \label{fig:4_1_7}
     \end{subfigure}
        \caption{Parameter initialization based on the optimal quadrature rule for a space with one fewer element, while keeping the degree and continuity. Red circles correspond to the points that are scaled for the initialization, while the blue circles correspond to the rest of the points in the rule. The area of each circle is proportional to its weight $w_i$.}
        \label{fig:uniform_initialization}
\end{figure}

We consider a partition $U$ to be uniform if the distance between any two consecutive non-repeating knots $u_{j-1}$ and $u_j$, ${j=1,\ldots,\ne}$, is constant. In these cases, due to symmetry in the function space around the center, we can impose two additional constraints to our optimal quadrature rule:
\begin{equation}
\begin{aligned}
 \{x_i\}^\frac{q}{2}_{i=1} &= \{1-x_j\}^\frac{q}{2}_{j=q} \,,\\[2pt]
 \{w_i\}^\frac{q}{2}_{i=1} &= \{w_j\}^\frac{q}{2}_{j=q} \,.
 \label{eq:symmetry_constraints}
\end{aligned}
\end{equation}
Thus, the number of parameters to fit is $\phi_x = \lfloor \frac{q}{2} \rfloor$ points and $\phi_w = \lceil \frac{q}{2} \rceil - 1$ weights. Then, given a space $\Sr_k^{d}\big(\{u_j\}^{\ne}_{j=0}\big)$, for ${k>0}$, the parameter initialization will be done as follows:
\begin{align}
 \{x_i\}^{\phi_x}_{i=1} &= \widetilde{x}_i\cdot \frac{\ne-1}{\ne} \,,\\[2pt]
 \{w_i\}^{\phi_w}_{i=1} &= \widetilde{w}_i\cdot \frac{\ne-1}{\ne} \,,
\end{align}
where $\{\widetilde{x}_i\}^{\widetilde{q}}_{i=1}$ and $\{\widetilde{w}_i\}^{\widetilde{q}}_{i=1}$ are the optimal quadrature points and weights for the space $\Sr_k^{d}\big(\{u_j\}^{\ne-1}_{j=0}\big)$. That is, the initial values when searching for a solution in a partition with $\ne$ elements result from a linear scaling of the optimal solution for a partition with $\ne-1$ elements. The base case is $\ne=2$, in which the initial values for $x_i$ are uniformly distributed in the interval $[0, 0.5]$, and all the weights are $w_i=\frac{1}{q}$. This procedure is illustrated in \autoref{fig:uniform_initialization}. Even if the number of optimal quadrature points increases with the number of elements (for fixed $d$ and $k$), it holds that $\widetilde{q} \geq \frac{q}{2}$, and therefore this procedure can always be applied.

For zero-continuity spaces (i.e., ${k=0}$), a better initialization is achieved as follows:
\begin{align}
 x_1 = \frac{1}{2q}\,, \qquad \{x_i\}^{\phi_x}_{i=2} &= x_{i-1} + \frac{1}{q} \,,\\[2pt]
                            \{w_i\}^{\phi_w}_{i=1} &= \frac{1}{q} \,.
\end{align}
The main advantage in this case is that the initialization does not depend on a quadrature rule for a different space, and therefore we can directly run the optimization for an arbitrary number of elements without having to solve the problem for fewer elements.

\subsubsection{Non-uniform elements}
\label{sub:non-uniform_init}

\begin{figure}[!h]
     \centering
     \begin{subfigure}[b]{0.48\textwidth}
         \centering
         \includegraphics[width=\textwidth]{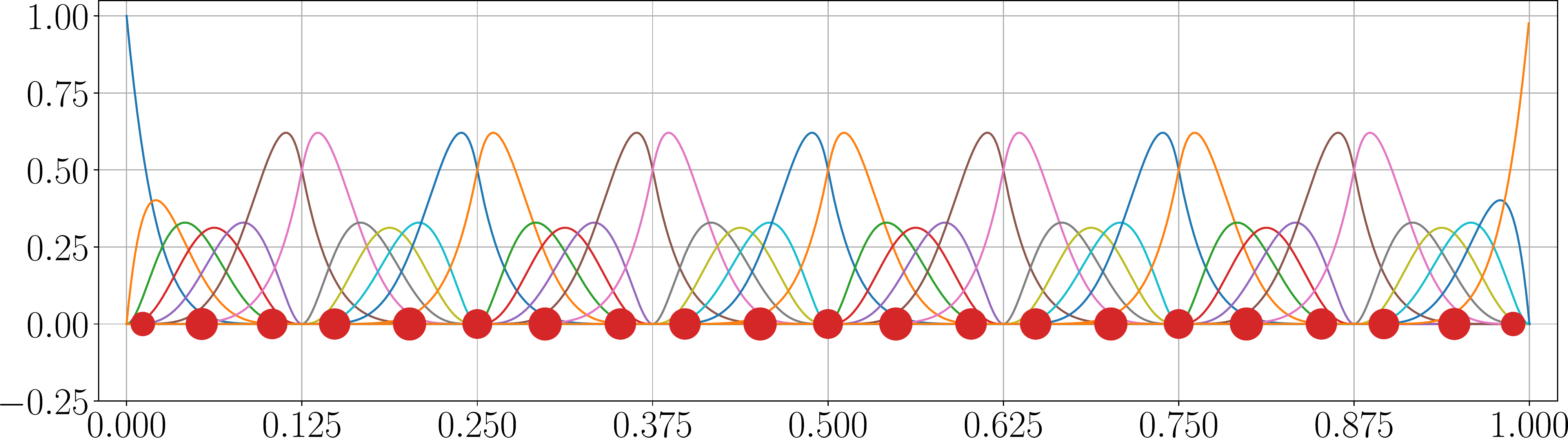}
         \caption{Optimal quadrature points for a uniform partition with $d=6, k=1, \ne=8$.}
         \label{fig:4_1_8}
     \end{subfigure}
     \hfill
     \begin{subfigure}[b]{0.48\textwidth}
         \centering
         \includegraphics[width=\textwidth]{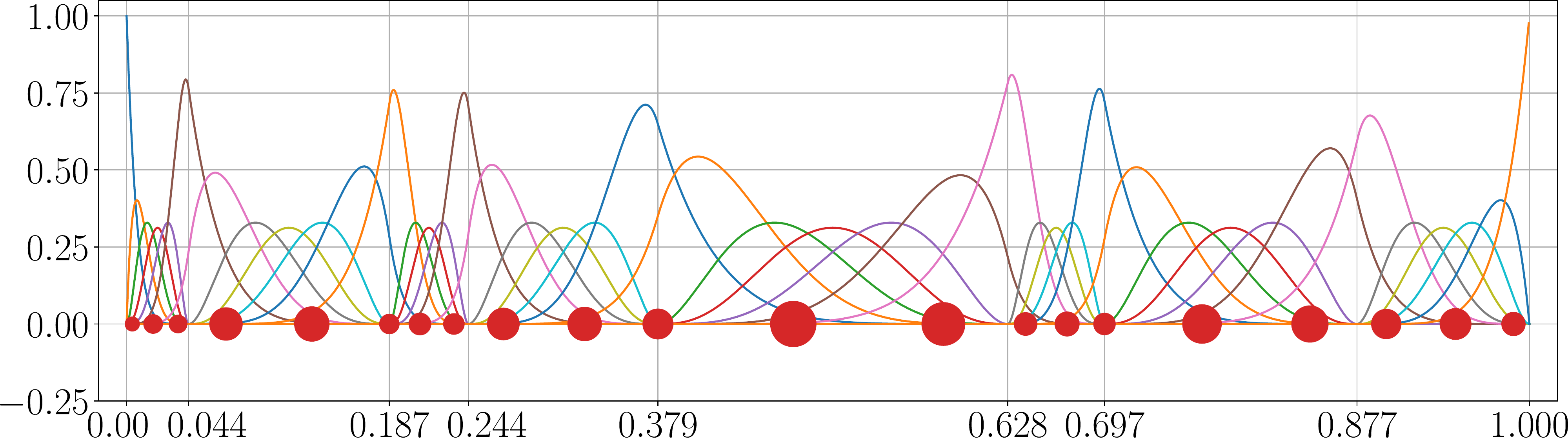}
         \caption{Initial parameter values for a random non-uniform partition with $d=6, k=1, \ne=8$.}
         \label{fig:4_1_8_nunif}
     \end{subfigure}
        \caption{Parameter initialization for non-uniform partitions based on the optimal quadrature rule for a space with a uniform partition and the same degree, continuity, and number of elements. The basis functions are shown as a reference, and the circles correspond to the points and weights of the quadrature rule. The area of each circle is proportional to its weight $w_i$.}
        \label{fig:nonuniform_initialization}
\end{figure}

When our element domain forms a non-uniform partition, we cannot exploit the constraints in~\eqref{eq:symmetry_constraints} to reduce the number of optimization parameters because the symmetry of the function space around the center no longer holds.
In such occasions, we initialize the values of $x_i$ and $w_i$ from the optimal rule for a space with a uniform partition and the same $d$, $k$ and $\ne$. Thus, given a space $\Sr_k^{d}\big(\{u_j\}^{\ne}_{j=0}\big)$ with a non-uniform partition (i.e., $u_j-u_{j-1}$ is not constant), and the corresponding space $\Sr_k^{d}\big(\{\widetilde{u}_j\}^{\ne}_{j=0}\big)$ with a uniform partition and its optimal quadrature rule as $\{\widetilde{x}_i\}^q_{i=1}$, $\{\widetilde{w}_i\}^q_{i=1}$, we calculate the initial parameter values for every $\widetilde{x}_i \in [\widetilde{u}_j, \widetilde{u}_{j+1}]$ as follows:
\begin{align}
  x_i &= (\widetilde{x}_i-\widetilde{u}_j)\cdot\frac{u_{j+1}-u_j}{\widetilde{u}_{j+1}-\widetilde{u}_j} + u_j \,,\\[2pt]
  w_i &= \widetilde{w}_i\cdot\frac{u_{j+1}-u_j}{\widetilde{u}_{j+1}-\widetilde{u}_j}\,.
\end{align}
Basically, we calculate a \textit{length scale} given by $\frac{u_{j+1}-u_j}{\widetilde{u}_{j+1}-\widetilde{u}_j}$ for each element with respect to the uniform partition, and then linearly apply this scale to the optimal points and weights belonging to that element. \autoref{fig:nonuniform_initialization} visually shows an example for a space with a random non-uniform partition of 8 elements.

\subsection{Optimization loop}
\label{sub:optimization}

Once the initial values are set, we start a regular gradient-descent search using the Yogi optimizer~\cite{NEURIPS18_adaptive}. The only parameter of the optimizer that has to be tuned is the learning rate. We notice that this tuning should be driven by the complexity of the problem represented by $q$, since the larger the number of parameters, the more sensitive the loss function is to small numeric variations, affecting the smoothness of the search space. We also observed that with just an inverse linear scaling, the search occasionally stalls by jumping around the solution, and thus we opt for an inverse log-linear scaling by setting the learning rate as $\frac{10^{-2}}{q\cdot \log(q)}$.

On the other hand, since an exact solution is difficult to find in practice through a gradient-descent approach (due to floating-point arithmetic), we consider a quadrature rule to be optimal if ${L(X, W) < 10^{-20}}$. The search stops if either ${L(X, W) < 10^{-25}}$ or we reach the limit of 10,000 gradient-descent iterations (a.k.a epochs).


\section{Numerical results}
\label{sec:Results}

\subsection{Implementation details}

The proposed method has been fully implemented in the Python programming language, relying heavily on the JAX~\cite{jax2018github} and Optax~\cite{deepmind2020jax} libraries. All the experiments have been performed in a laptop PC with an Intel(R) Core(TM) i7-1250U CPU, 32~GB of RAM memory and a Linux kernel v5.19. The version of Python in which we ran the experiments was 3.10.

With the aim of supporting reproducible research, the full source code required to reproduce the experimental results described below is provided as open-source software. We also provide a table with the optimal quadrature rules in double precision for spaces up to $d=16$ and arbitrary $k$ with uniform partitions of up to ${\ne=50}$ elements. All this may be downloaded from \url{https://gitlab.bcamath.org/tteijeiro/spline-integration}

\subsection{General experiments}
\label{sub:GeneralExperiments}

\begin{figure}[!h]
    \centering
    \includegraphics[width=\textwidth]{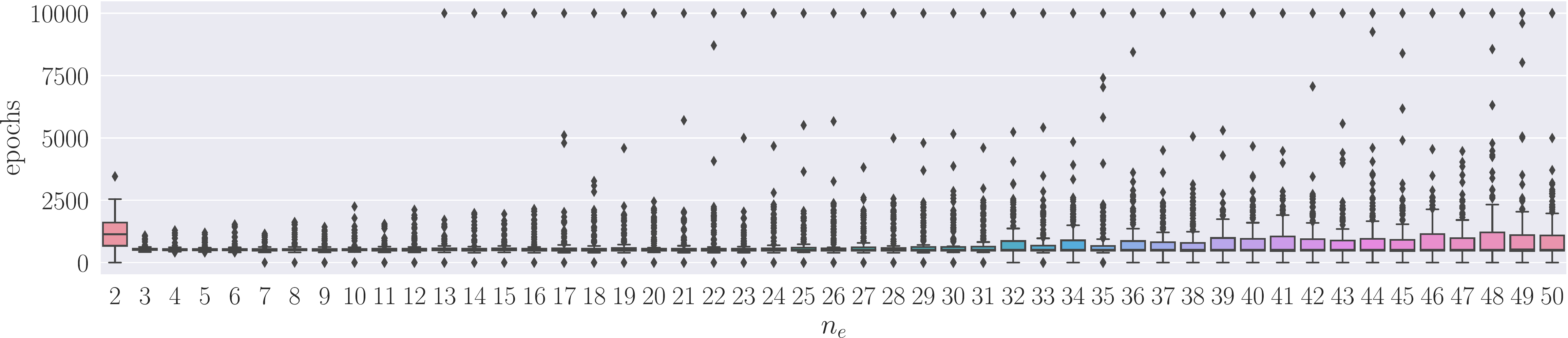}
    \caption{Boxplot with the distributions of the number of epochs required to converge to the solution for different number of elements. Each box spans trough the quartiles of the distribution, and the whiskers span up to 1.5 times the interquartile range. The rest of points are individually plotted as outliers. The black horizontal bar within each box corresponds to the median.}
    \label{fig:nepoc_distribution}
\end{figure}

The first experiments aim at evaluating the generality of the proposed method, both with uniform and non-uniform partitions. For this, we performed an exhaustive search of optimal quadrature rules for spline spaces with all combinations of ${k< d \leq 16}$ and uniform partitions with ${2 \leq \ne \leq 50}$, corresponding to 6664 cases. For all cases, we succeed in obtaining optimality, so that ${L(X,W) < 10^{-20}}$. \autoref{fig:nepoc_distribution} shows the distribution of the number of epochs required to converge to the optimal rule based on the number of elements. Even if the difficulty clearly increases with the number of elements, the scalability shown by the method is really positive, with a global median of 522 required epochs per rule, and with 90\% of the rules requiring less than 1700 epochs to meet the stopping criteria of ${L(X, W) < 10^{-25}}$. It is also interesting to see the case of ${\ne=2}$ in which the recursive parameter initialization is not applied and the number of required epochs is generally higher, thus indicating the advantage of the proposed initialization technique.

\begin{figure}[h]
     \centering
     \begin{subfigure}[b]{0.48\textwidth}
         \centering
         \includegraphics[width=\textwidth]{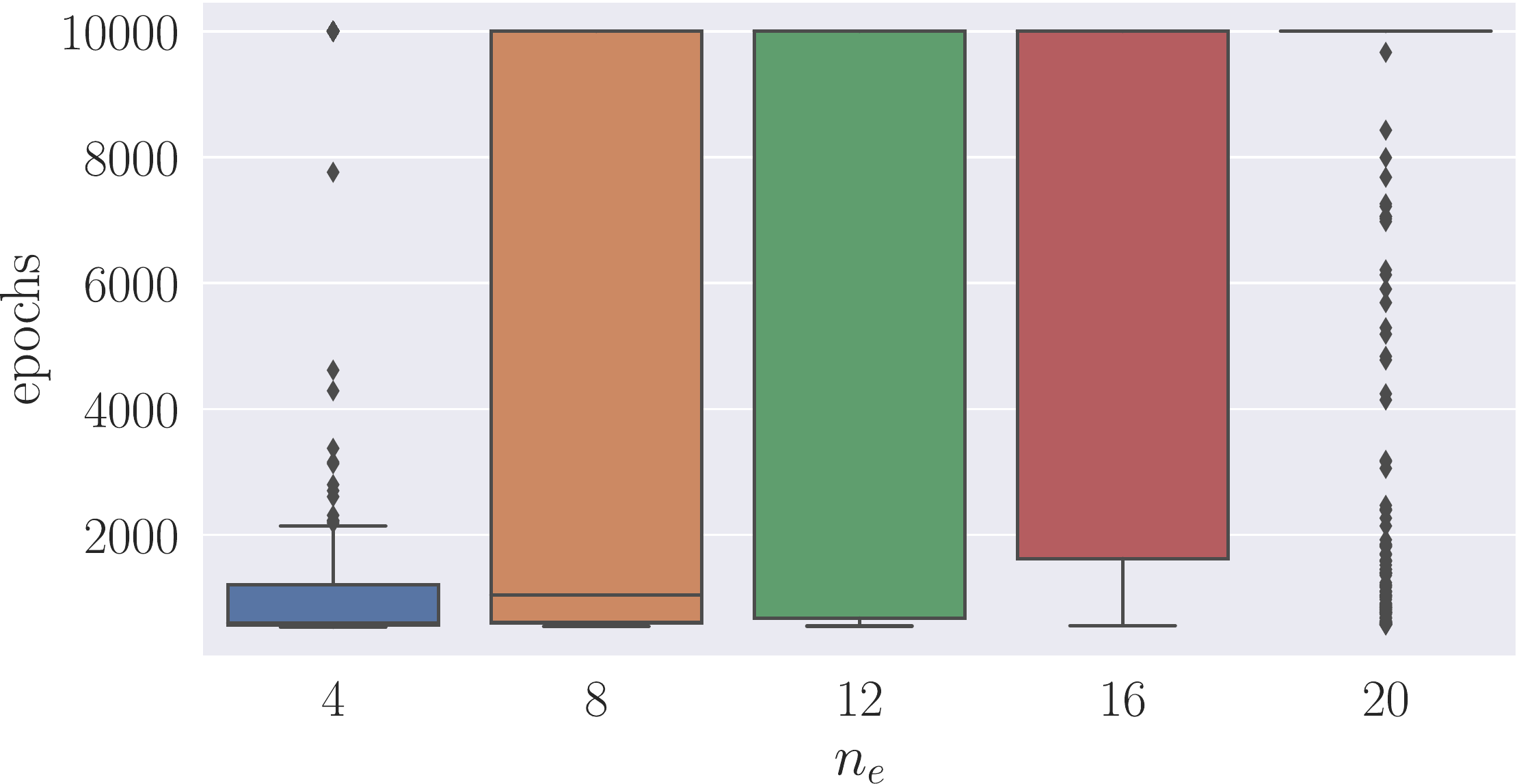}
         \caption{Boxplot with the distributions of the number of optimization epochs run for different number of elements in non-uniform partitions.}
         \label{fig:nepoc_distribution_nunif}
     \end{subfigure}
     \hfill
     \begin{subfigure}[b]{0.48\textwidth}
         \centering
         \includegraphics[width=\textwidth]{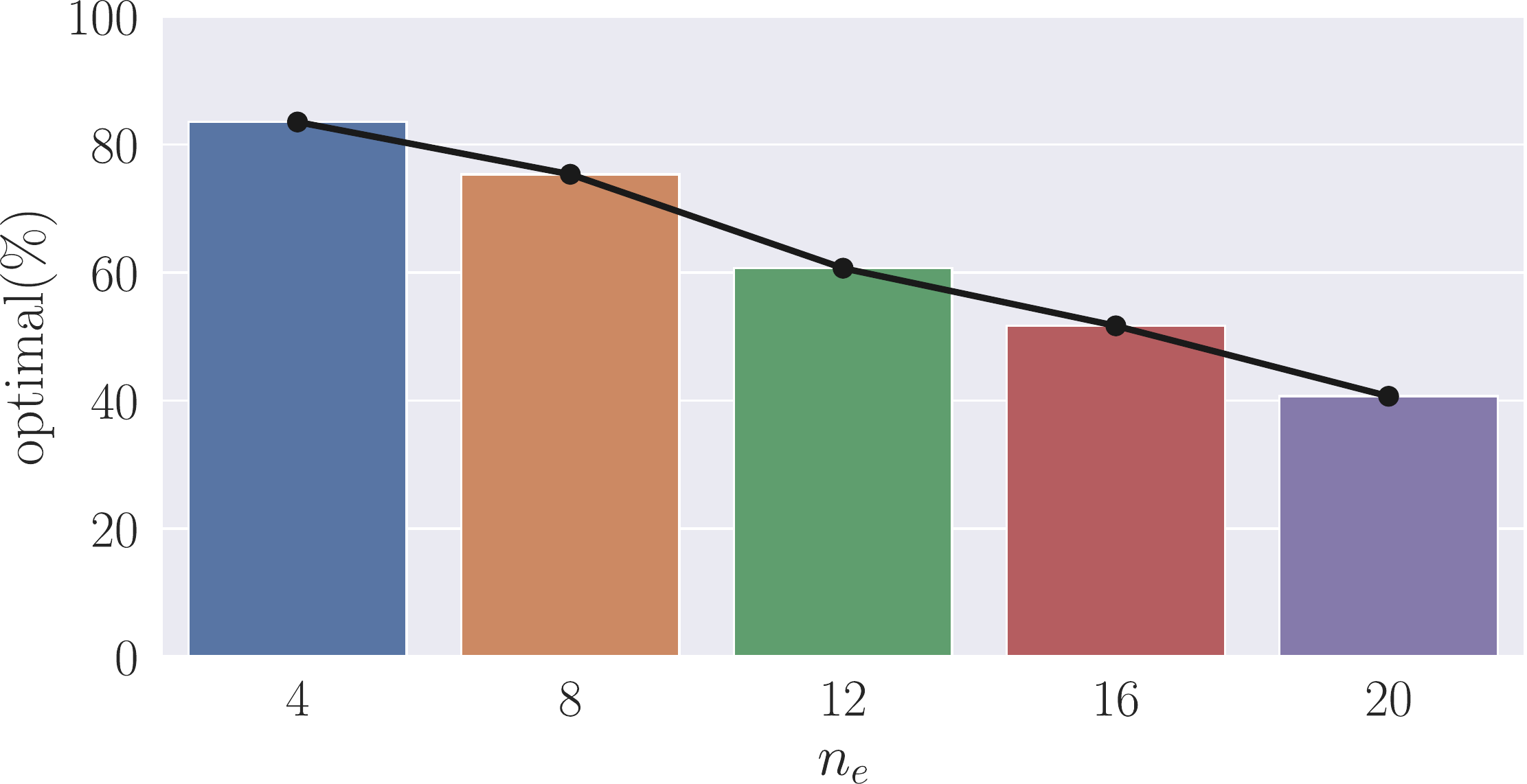}
         \caption{Percentage of non-uniform partitions for which an optimal quadrature rule was found, according to the number of elements.}
         \label{fig:nunif_success_rate}
     \end{subfigure}
        \caption{Results on finding optimal rules for non-uniform partitions with different number of elements.}
        \label{fig:nonuniform_results}
\end{figure}

To test the method on non-uniform partitions, we randomly generate 3000 function spaces with the following values of $(d, k)$: $(4, 0), (5, 0), (6, 1), (7, 1), (8, 2)~\text{and}~(9, 2)$, and with a variable number of elements ${\ne=4, 8, 12, 16, 20}$. For each combination of $d, k~\text{and}~n_e$, we generate 100 independent partitions by randomly sampling and sorting ${n_e-1}$ values from a uniform distribution in the interval $[0, 1]$. Then, for each of these 3000 function spaces, we initialize and optimize the parameters $x_i$ and $w_i$ according to the method described in \az{Sections}~\ref{sub:non-uniform_init} and~\ref{sub:optimization}.
Unlike with uniform partitions, in this case the optimization loop did not converge to an optimal quadrature rule in 100\% of the cases. \autoref{fig:nonuniform_results} shows the results for different number of elements. \autoref{fig:nepoc_distribution_nunif} depicts the distributions in the number of epochs, which are in general much higher than uniform partitions. Additionally, more than 50\% of the spaces reach the limit of 10000 epochs. \autoref{fig:nunif_success_rate} shows the percentage of partitions for which the method was able to find an optimal rule, that is~${L(X,W)<10^{-20}}$, according to the number of elements. The reduction in the success rate is fundamentally linear, and thus the method cannot be considered of general applicability for ${\ne > 8}$. In future work we will focus on characterizing additional constraints that must be satisfied by the partition knots to improve the convergence.

\subsection{Computational requirements}

\begin{figure}[!h]
    \begin{subfigure}{0.33\textwidth}\hspace{-0.3cm}
        \includegraphics{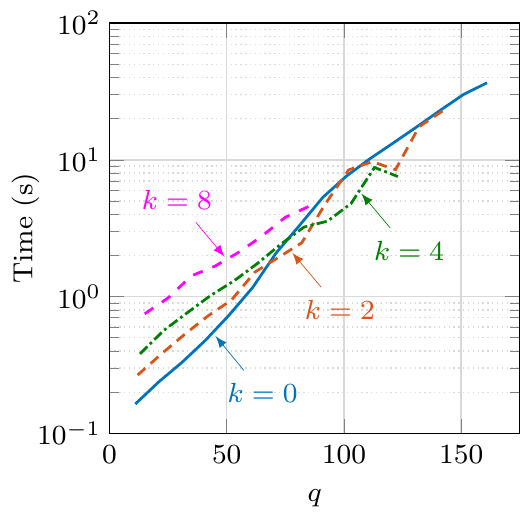}
        \caption{Increasing $q$ by increasing $d$.}
        \label{fig:time_deg}
    \end{subfigure}
    \begin{subfigure}{0.33\textwidth}
        \includegraphics{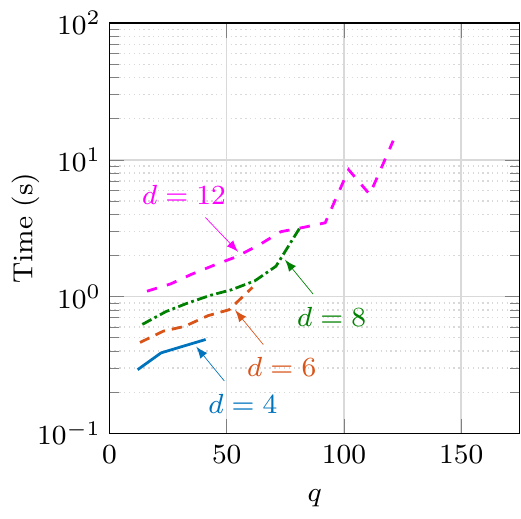}
        \caption{Increasing $q$ by decreasing $k$.}
        \label{fig:time_cont}
    \end{subfigure}
    \begin{subfigure}{0.33\textwidth}\hspace{0.1cm}
        \includegraphics{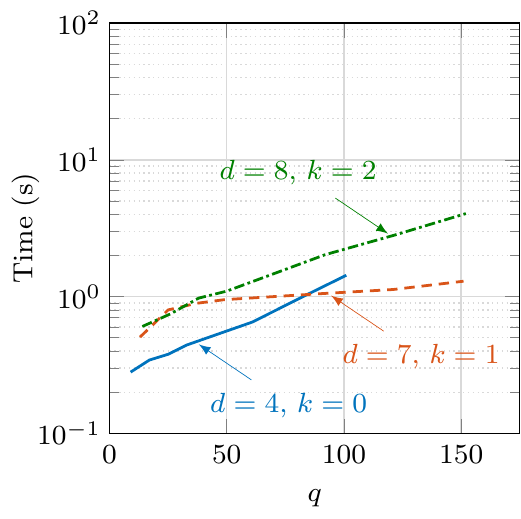}
        \caption{Increasing $q$ by increasing $n_e$.}
        \label{fig:time_ne}
    \end{subfigure}
    \caption{Variation of the computational time required to calculate a quadrature rule on uniform partitions according to the number of parameters $q$. Since $q$ is determined by $d$, $k$, and $n_e$, we show the influence on each of these variables separately.}
    \label{fig:timing_results}
\end{figure}

To evaluate the computational demands of our method, we assess the time required to find a quadrature rule with respect to the complexity of the rule, measured as the number of quadrature points to optimize~($q$) and considering only uniform partitions. However, $q$ is determined by the spline space, and there might be two different spaces leading to the same $q$. Thus, we will separately assess the influence of $c$, $d$, and $n_e$ to get a better overview of how the method scales in practice.
\autoref{fig:timing_results} shows three graphs of the time required to calculate rules for increasing values of $q$. In \autoref{fig:time_deg}, the number of points $q$ increases by increasing the maximum degree of the space $d$, while keeping a fixed number of elements ${n_e=20}$ and a fixed continuity~$c$ for each plotted line. On the other hand, in \autoref{fig:time_cont} the value of $q$ increases by decreasing $c$, keeping the same ${n_e=20}$ and a fixed $d$ for each line. Finally, \autoref{fig:time_ne} keeps a fixed degree and continuity for each line, and increases $n_e$ between 4 and 50 to increase $q$. All charts are in logarithmic scale on the time axis.

We can see that, in general, the time requirements grow exponentially with the number of parameters. Still, this growth is much faster when the reason is the complexity of the function space (higher degree or lower continuity) than the case where the reason is an increased number of elements. Indeed, we can observe that for a space with $d=7$ and $c=1$, the increase is basically linear, and this is due to the fact that, in this space, the initialization method described in \autoref{sub:initialization_uniform} already gives an optimal rule, not requiring any gradient descent iterations. From a practical point of view, we can see that the quadrature rule for a realistic space with $d=7$, $c=1$ and $n_e=50$ can be calculated in approximately 1 second, while an overly complex space with $d=16$, $c=0$ and $n_e=20$ can take up to 36 seconds. We want to reiterate that these experiments were performed in a common laptop CPU, so these times could be significantly reduced by using a more powerful hardware.

\subsection{Optimal vs. element-wise Gaussian rules}
\label{sub:OptimalvsEWG}

\begin{figure}[!h]\centering
	\begin{subfigure}{0.49\textwidth}\centering
        \includegraphics{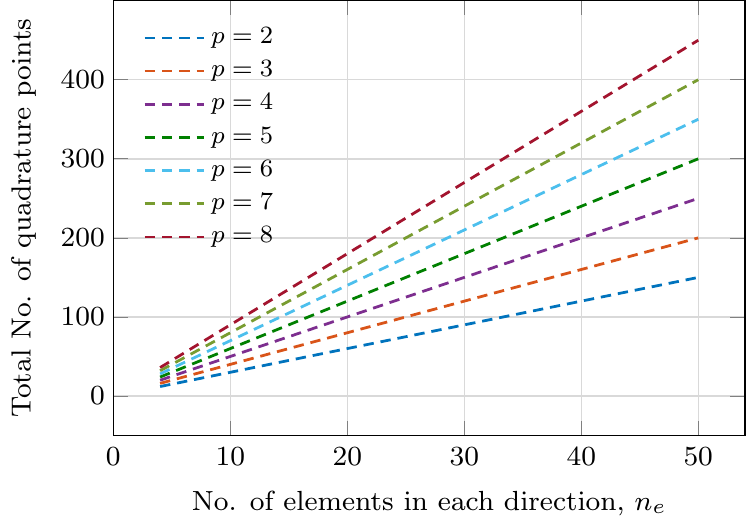}
        \caption{EWG rule in 1D.}
	\end{subfigure}
	\begin{subfigure}{0.49\textwidth}\centering
        \includegraphics{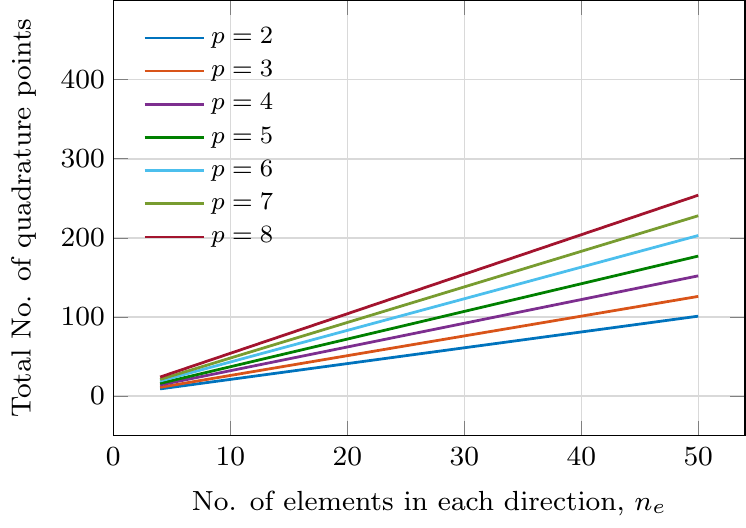}
        \caption{Optimal rule in 1D.}
	\end{subfigure}
	\begin{subfigure}{0.49\textwidth}\centering
        \includegraphics{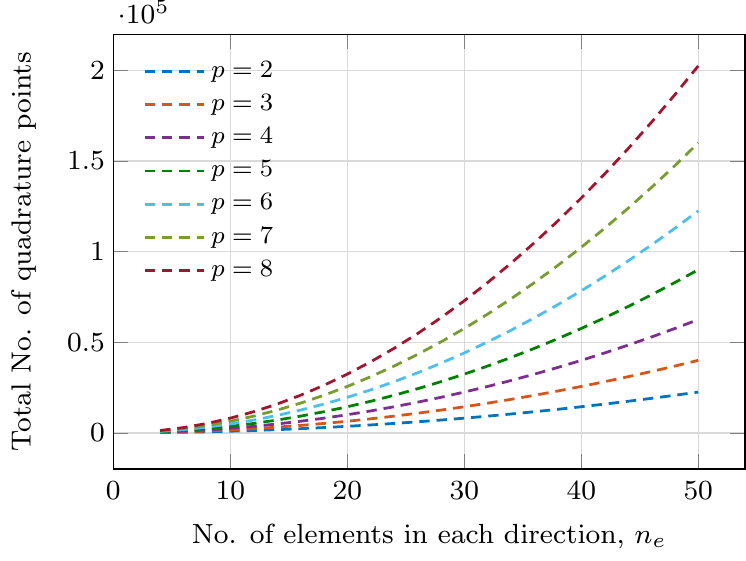}
        \caption{EWG rule in 2D.}
	\end{subfigure}
	\begin{subfigure}{0.49\textwidth}\centering
        \includegraphics{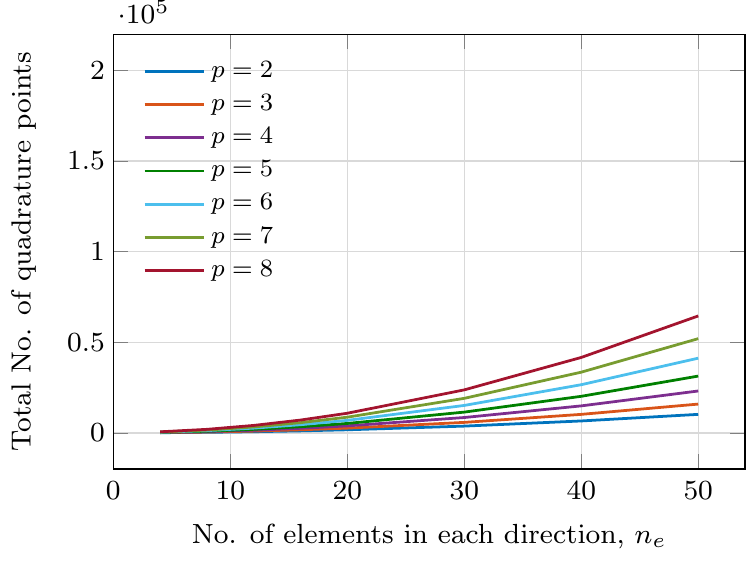}
        \caption{Optimal rule in 2D.}
	\end{subfigure}
	\begin{subfigure}{0.49\textwidth}\centering
        \includegraphics{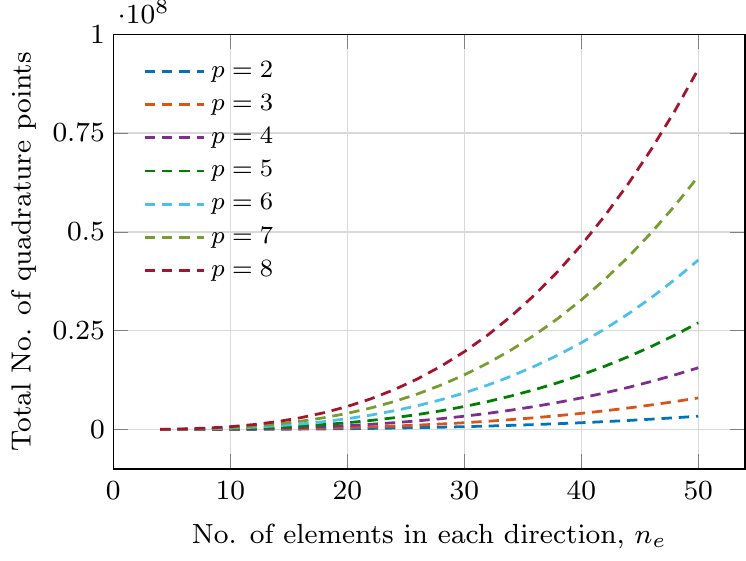}
        \caption{EWG rule in 3D.}
	\end{subfigure}
	\begin{subfigure}{0.49\textwidth}\centering
        \includegraphics{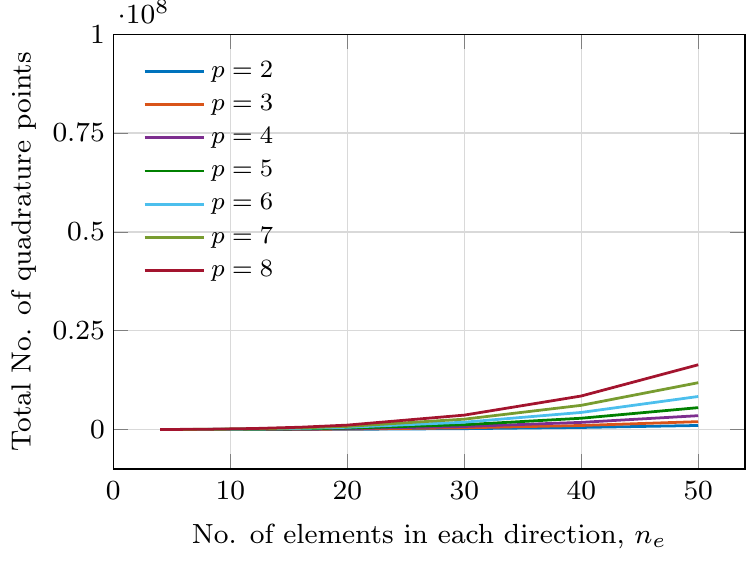}
        \caption{Optimal rule in 3D.}
	\end{subfigure}
    \caption{Comparison of the optimal and EWG quadrature rules in terms of the total number of required quadrature points in different dimensions.}
    \label{fig:IntegrationCost}
\end{figure}

Herein, we compare the computational cost of numerical integration when using the optimal and classical EWG quadrature rules. For this purpose, we consider the total number of required quadrature points in 1D, 2D, and 3D parameter spaces when constructing the mass and stiffness matrices in the sense of IGA discretization.
For 2D and 3D spaces, when the coefficients of the desired PDE are constant along different spatial directions, one may construct the system matrices in 1D and then obtain the matrices in higher dimensions through the tensor product of 1D matrices (see, e.g.,~\cite{Hashemian2021,Puzyrev2017}). In this case, the cost of numerical integration in higher dimensions is equal to that of 1D spaces. However, when the coefficients are varying through the domain in different directions, which could be the case in problems arising in electromagnetics and geoscience (see, e.g.,~\cite{Garcia2019,Hashemian20212}), we construct the system matrices by integrating higher-dimensional basis functions. Thus, the required number of quadrature points (i.e., the integration cost) grows by the power of the dimension space and the advantage of using the proposed optimal quadrature rules is more tangible. \autoref{tab:IntegrationCost} shows the specific savings for some relevant spaces in different dimensions when using functions in $\Sr^{2p}_{p-2}(U)$ as the integrand in the sense of \autoref{rem:Spaces}. The improvements for all tested values of $p$ and $n_e$ are illustrated in \autoref{fig:IntegrationCost}.

\begin{table}[!h]\centering
	\caption{Comparison of the required number of quadrature points for optimal and element-wise Gaussian rules when constructing stiffness and mass matrices using the IGA discretization. We report results for some selected polynomial degrees and number of elements in 1D, 2D and 3D spaces.}
	\label{tab:IntegrationCost}
    \small
	\begin{tabular}{@{}cccrrc@{}}
		\toprule
		\multirow{2}{*}{Dimensionality}                 & Degree             & Number of elements & \multicolumn{2}{c}{Number of quadrature points} & \multirow{2}{*}{\shortstack{Relative savings when\\ using optimal rules}}\\[0.7em]
		      & $p$                & $n_e$              & EWG rule              & Optimal rule            &   \\ \midrule
		\multirow{6}{*}{1D} & \multirow{2}{*}{2} & 20                 & 60                    & 41                      & 31.7\%                \\
		&                    & 50                 & 150                   & 101                     & 32.7\%                \\ \cline{2-6} 
		& \multirow{2}{*}{4} & 20                 & 100                   & 62                      & 38.0\%                \\
		&                    & 50                 & 250                   & 152                     & 39.2\%                \\ \cline{2-6} 
		& \multirow{2}{*}{8} & 20                 & 180                   & 104                     & 42.2\%                \\
		&                    & 50                 & 450                   & 254                     & 43.6\%                \\ \midrule
		\multirow{6}{*}{2D} & \multirow{2}{*}{2} & 20                 & 3,600                 & 1,681                   & 53.3\%                \\
		&                    & 50                 & 22,500                & 10,201                  & 54.7\%                \\ \cline{2-6} 
		& \multirow{2}{*}{4} & 20                 & 10,000                & 3,844                   & 61.6\%                \\
		&                    & 50                 & 62,500                & 23,104                  & 63.0\%                \\ \cline{2-6} 
		& \multirow{2}{*}{8} & 20                 & 32,400                & 10,816                  & 66.7\%                \\
		&                    & 50                 & 202,500               & 64,516                  & 68.1\%                \\ \midrule
		\multirow{6}{*}{3D} & \multirow{2}{*}{2} & 20                 & 216,000               & 68,921                  & 68.1\%                \\
		&                    & 50                 & 3,375,000             & 1,030,301               & 69.5\%                \\ \cline{2-6} 
		& \multirow{2}{*}{4} & 20                 & 1,000,000             & 238,328                 & 76.2\%                \\
		&                    & 50                 & 15,625,000            & 3,511,808               & 77.5\%                \\ \cline{2-6} 
		& \multirow{2}{*}{8} & 20                 & 5,832,000             & 1,124,864               & 80.7\%                \\
		&                    & 50                 & 91,125,000            & 16,387,064              & 82.0\%                \\ 
		\bottomrule
	\end{tabular}
\end{table}

\section{Case studies}
\label{sec:CaseStudies}

To study the accuracy of the proposed optimal quadrature rules on spline spaces, we consider two case studies, namely, the eigenproblem of the Laplace operator and the eigenfrequency analysis of a curved beam.

\subsection{Eigenvalue problem of the Laplace operator}

We consider the eigenproblem of the Laplace operator in 1D. Given ${\Om:[0,1]\subset\R}$ as our computational domain with boundary $\partial\Om$, we write
\begin{align}
    \begin{cases}
        {\rm Find}~\lm\in\R^+~{\rm and}~\ur:\Om\rightarrow\R,~{\rm such~that}\\[2pt]
        \qquad
        \begin{aligned}
            \Delta\ur + \lm\ur &= 0\,,\quad &&{\rm in}~\Om\,,\\
                           \ur &= 0\,,\quad &&{\rm on}~\partial\Om\,.
        \end{aligned}
    \end{cases}
    \label{eq:Laplace}
\end{align}
Let us consider ${H^1(\Om)}$ gradient-conforming functional space. We define the Sobolev space ${H^1_0(\Om)}$ of functions in ${H^1(\Om)}$ that are vanishing on the boundary.
Taking ${\vr\in H^1_0(\Om)}$ as an arbitrary test function, we build the weak form of~\eqref{eq:Laplace} as follows:
Find ${\lm\in\R^+}$ and ${\ur\in H^1_0(\Om)}$, such that for all ${\vr\in H^1_0(\Om)}$,
\begin{align}
	a\.(\vr,\ur)=\lm\.b\.(\vr,\ur)\,,
	\label{eq:LaplaceWeak}
\end{align}
where
\begin{align}
	a\.(\vr,\ur) &:= \int_{\Om} \nabla\vr\cdot\nabla\ur \,d\Om\,,\\
	b\.(\vr,\ur) &:= \int_{\Om} \vr\ur \,d\Om\,.
\end{align}
Then, by employing an IGA discretization spanned by $\Sr^p_c(U)$ over the knot sequence ${U=\{u_j\}_{j=0}^\ne}$ with bases ${\{v_i\}_{i=0}^n}$\., we define the discrete space ${\Vr^h:=\Sr^p_c(U)}$. Thus, we obtain the discrete eigenproblem:
Find ${\lm^h\in\R^+}$ and ${\ur^h\in \Vr^h\subset H^1_0(\Om)}$, such that for all ${\vr^h\in \Vr^h\subset H^1_0(\Om)}$,
\begin{align}
	a\.\big(\vr^h,\ur^h\big)=\lm^h\.b\.\big(\vr^h,\ur^h\big)\,,
	\label{eq:LaplaceWeakDiscrete}
\end{align}
where $\lm^h$ and $\ur^h$ refer to approximated eigenpairs.
Equation~\eqref{eq:LaplaceWeakDiscrete} leads to the following generalized eigenproblem in the matrix form:
\begin{align}
    K\.\Ur^h=\lm^h M\.\Ur^h\,,
    \label{eq:EigProblem}
\end{align}
where $\Ur^h$ denotes the eigenvectors. We obtain the stiffness and mass matrices, $K$ and $M$, respectively, using integrals~\eqref{eq:Stiffness} and~\eqref{eq:Mass}.
In the following, we consider the eigensolution of~\eqref{eq:EigProblem} when constructing the system matrices using both the proposed optimal and the classical element-wise Gaussian quadrature rules.
For the former, we create our rules based on the spline space ${\Sr^{2p}_{p-2}(U)}$ as the integrand (c.f. \autoref{rem:Spaces}).
In all tests, we show that the optimal quadrature rules have the same order of exactness as the classical EWG rules while constructing matrices with lower computational efforts (see \autoref{sub:OptimalvsEWG}).

\subsubsection{Uniform domain}

\begin{figure}[!h]\hspace{-0.4cm}
	\includegraphics{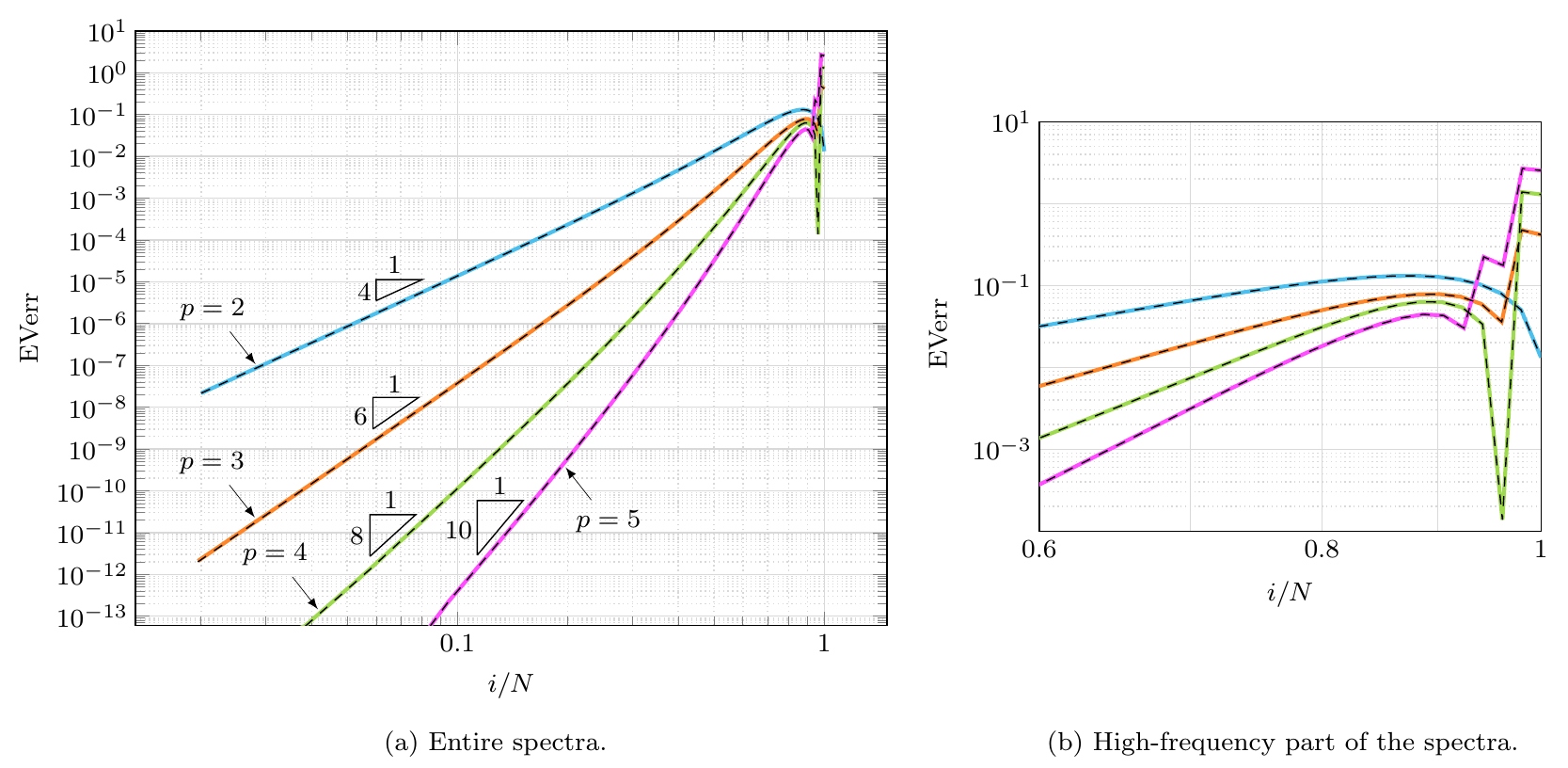}
	\caption{Eigenvalue error against the normalized mode number $i/\N$ (in logarithmic scale) when solving the eigenproblem of the Laplace operator discretized by ${\ne = 50}$ uniform elements. Solid lines refer to eigenvalue errors when using the proposed optimal rules. Dashed lines (identically superimposed over the solid lines) represent eigenvalue errors when using the classical EWG rules. The spectral approximation obtained by both approaches are identical (even in the high-frequency region where the outlier effect deteriorates the quality of the approximation) and in accordance with the theoretical convergence rates.}
	\label{fig:EVErr1D}	
\end{figure}

We first consider IGA discretizations with ${\ne=50}$ uniform elements and basis functions of polynomial degrees ${p=2,3,4,5}$. For every eigenpair, we introduce the eigenvalue error
\begin{align}
	{\rm EVerr} := \dfrac{\lm^h_i-\lm_i}{\lm_i}\,,
\end{align}
where ${\lm_i=i^{\.2}\pi^{\.2}}$, ${i=1,2,\ldots,\N}$, are the analytical eigenvalues, being $\N$ the total number of degrees of freedom (equal to ${n-2}$ in this 1D test case).
\autoref{fig:EVErr1D} shows, in logarithmic scale, the eigenvalue error against the normalized mode number $i/\N$ when solving~\eqref{eq:EigProblem} in a domain with ${\ne = 50}$ uniform elements. The results show that the spectral approximation obtained by both optimal and EWG rules are identical even in the high-frequency region where the outlier effect deteriorates the quality of the approximation. Since the domain consists of uniform elements, the theoretical convergence rate of $\Or(2p)$ is conserved (c.f.~\cite{Puzyrev2017}).

\subsubsection{Uniform domain with a fine mesh}

\begin{figure}[!h]\hspace{-0.2cm}
	\includegraphics{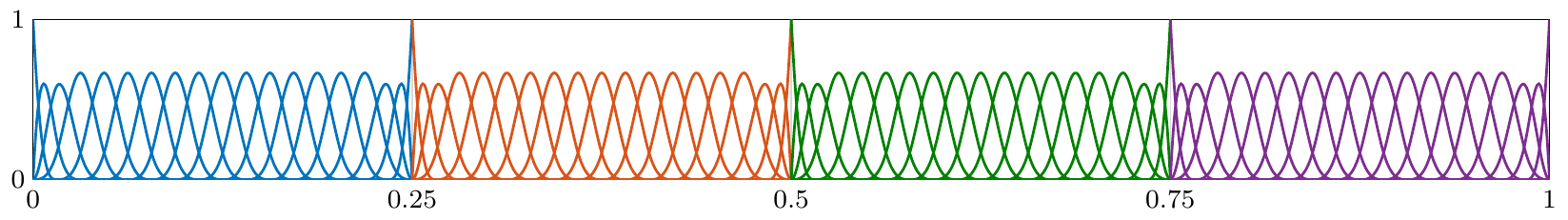}
	\caption{A domain with ${\ne=64}$ elements discretized by cubic basis functions. We subdivide the domain into four blocks of~16 elements and use the optimal quadrature rule generated for ${\Sr^6_1\big(\{u_j\}_{j=0}^{16}\big)}$ at each block in the sense of the rIGA discretization. The $C^0$ separator bases are shared between adjacent blocks.}
	\label{fig:rIGABasis}	
\end{figure}

In the next test, we consider finer domains (i.e., with higher number of elements) for which we do not construct any optimal rule.
Dealing with fine meshes is a common practice in FEA and IGA when we seek lower approximation errors.
In such occasions, we consider one optimal rule from those described in \autoref{sub:GeneralExperiments}. Then, we subdivide our domain into blocks of \textit{macroelements} in such a way that the number of elements at each block is equal to the number of elements for which we have the optimal rule.
This follows the terminology of the refined isogeometric analysis (rIGA)~\cite{Garcia2017}, where we introduce zero-continuity basis functions to split the domain into macroelements.
Since the system matrices at each block are constructed separately, we exploit the advantage of using the optimal rules for any given fine mesh with ${\ne>50}$ elements.
Furthermore, using $C^0$ separators in the mesh helps to reduce the interconnection between degrees of freedom, thus, reducing the solution costs when employing direct solvers for matrix factorization (c.f.~\cite{Garcia2017,Hashemian2021} for more details).
\autoref{fig:rIGABasis} depicts an example set of cubic basis functions in a domain with ${\ne=64}$ elements, that is subdivided into four blocks each containing 16 elements.

\begin{figure}[!h]\hspace{-0.4cm}
	\includegraphics{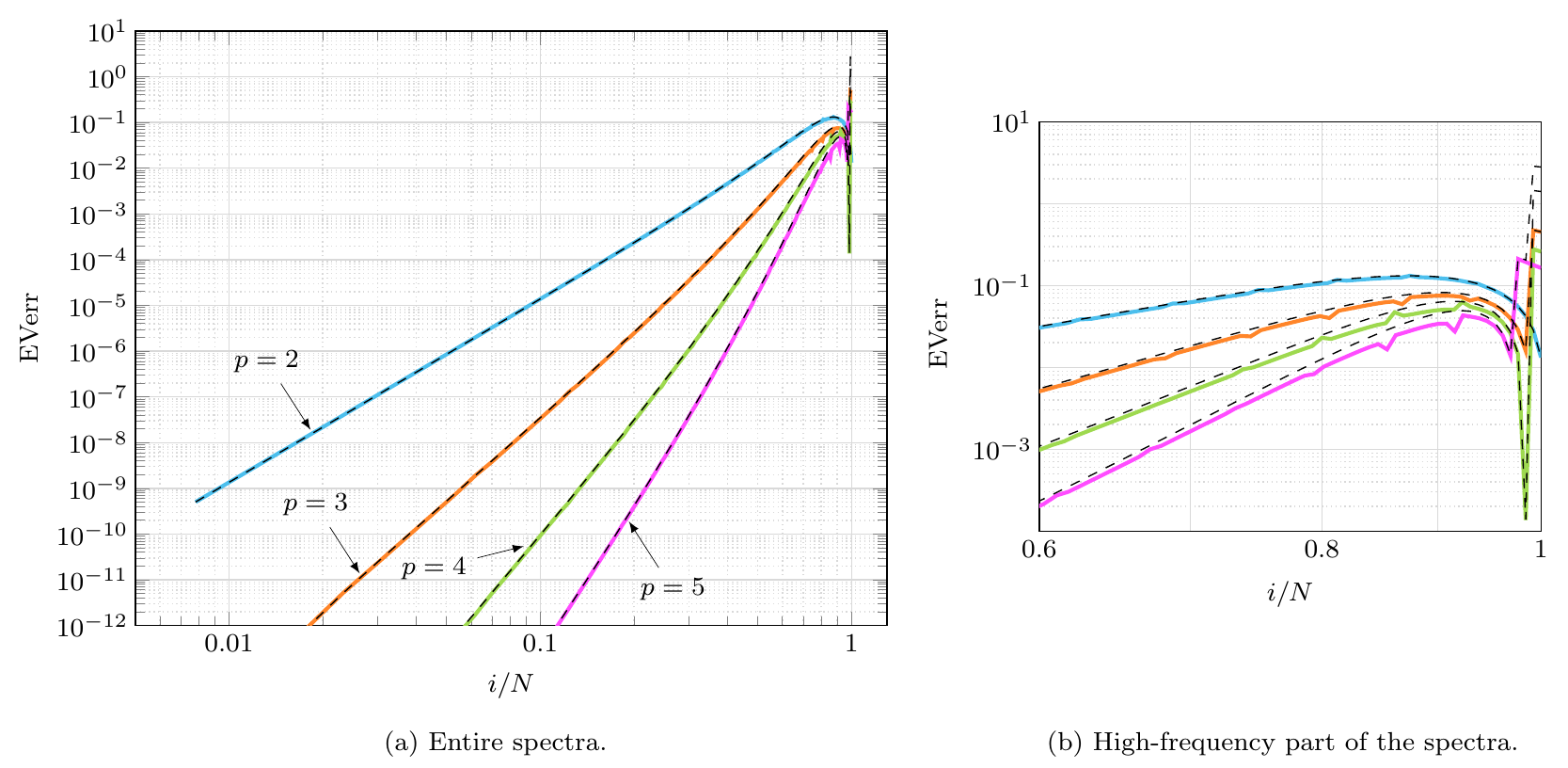}
	\caption{Eigenvalue error against the normalized mode number $i/\N$ for the solution of the eigenproblem of the Laplace operator discretized by ${\ne = 128}$ uniform elements. Solid and dashed lines correspond to the optimal and EWG rules, respectively. When integrating by the proposed optimal rules, we use the rIGA framework with blocks of 16 elements and consider ${N=N_{\rm IGA}}$ for mode number normalization. In this context, the high-frequency region shows slightly better spectral approximations, particularly for higher degrees.}
	\label{fig:EVErr1DrIGA}	
\end{figure}

To assess how the employment of the rIGA framework affects the results, we consider a uniformly-spaced 128-element domain and investigate the eigensolution accuracy of the eigenproblem of the Laplace operator. 
\autoref{fig:EVErr1DrIGA} shows the eigenvalue errors against the normalized mode number $i/\N$. We note that an rIGA-discretized system has more degrees of freedom compared to its IGA counterpart. This is because the continuity reduction at separators increases the dimension of the spline space. Thus, to have a true comparison with an IGA-discretized system constructed by EWG quadrature rules, we consider ${N=N_{\rm IGA}}$ for mode number normalization.
Herein, we consider blocks with 16 element and polynomial degrees ${p=2,3,4,5}$.
The accuracy of the eigensolution for the lower-frequency part of the spectra are identical when using optimal and EWG rules. However, the high-frequency region shows slightly better spectral approximations, particularly for higher degrees, when incorporating the optimal rules into the IGA framework. This improvement is mainly due to the effect that rIGA has on the spectral approximation (see, e.g.,~\cite{ Hashemian2021, Puzyrev2018} for more details).

\subsubsection{Non-uniform domain}

\begin{figure}[!h]\centering
	\includegraphics{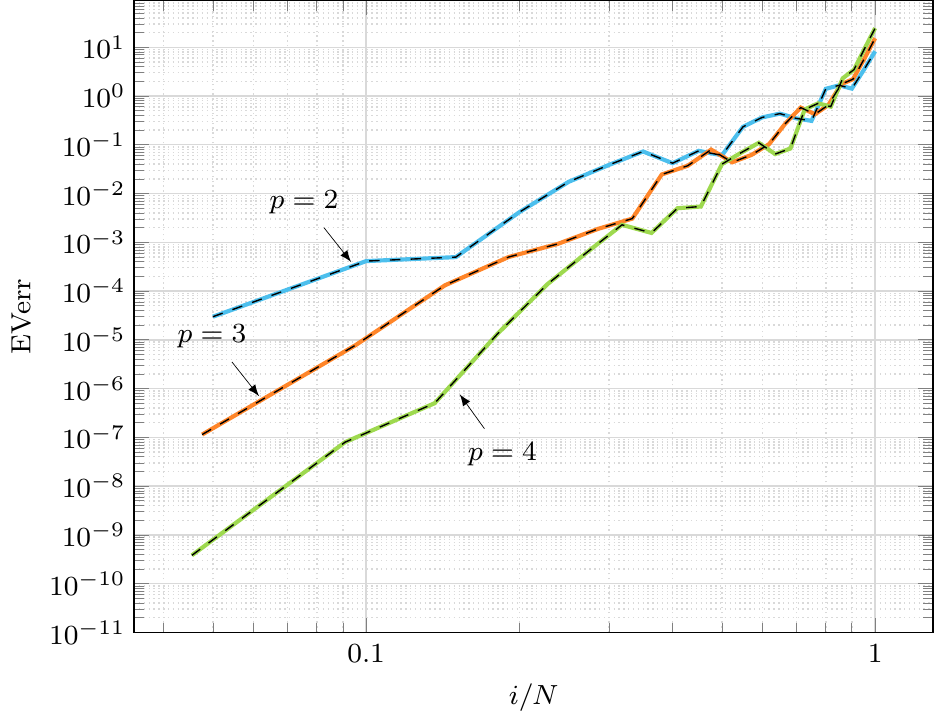}
	\caption{Eigenvalue error against the normalized mode number when solving the eigenproblem of the Laplace operator discretized by ${\ne = 20}$ non-uniform elements. The solid and dashed lines (correspond to the optimal and EWG rules, respectively) are identical, but deviate from the theoretical convergence rates due to the non-uniformity of the domain.}
	\label{fig:EVErr1Dnonuniform}	
\end{figure}

As a final test, we consider non-uniform elements in the IGA discretization of the eigenproblem~\eqref{eq:Laplace}.
Let us consider a domain with 20 non-uniform elements (for which our optimization is converged properly to optimal rules). In particular, we consider the solution over the following unevenly-spaced knot sequence:
\begin{multline}
U=[0.000,0.009,0.035,0.056,0.104,0.231,0.282,0.345,0.379,0.512,\\0.558,0.577,0.613,0.649,0.719,0.771,0.914,0.927,0.948,0.981,1.000]\..
\end{multline}
\autoref{fig:EVErr1Dnonuniform} shows the eigenvalue errors considering ${p=2,3,4}$ as the polynomial degree of basis functions.
Again, the accuracy of the spectral approximation is identical when using either of the optimal or EWG rules for system construction. However, it is clear that the eigenvalue error plots deviate from the theoretical convergence rates due to the non-uniformity of the domain.

\subsection{Eigenfrequency analysis of the Tschirnhausen beam, a curved geometry}

In the second case study, we investigate the performance of the proposed optimal quadrature rules when integrating over a curved geometry.
For this purpose, we consider the eigenfrequency analysis of the Tschirnhausen beam characterized by a planar freeform geometry with varying curvature.

\subsubsection{Geometry representation}

\begin{figure}[!h]\centering
	\includegraphics{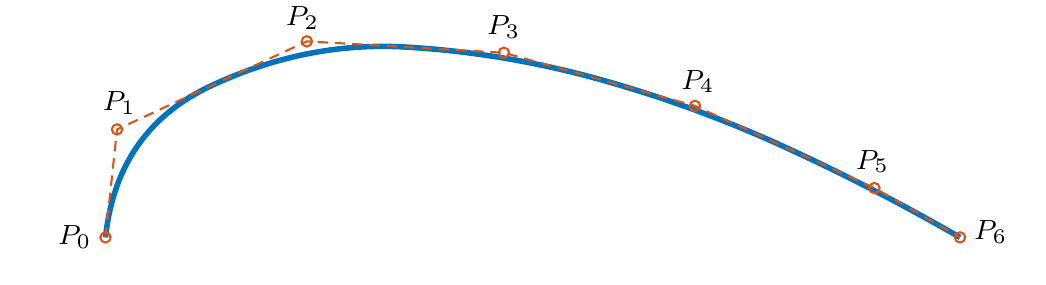}
	\caption{Geometry representation of the midline of the Tschirnhausen beam (blue) with seven control points (red dots).}
	\label{fig:Tschirn}	
\end{figure}

Given the parametric domain ${\hat{\Om}:[0,1]\subset\R}$, we represent the geometry of the beam midline by a B-spline curve ${C:\hat{\Om}\rightarrow\R^2}$ with ${n+1}$ control points ${P_i\in\R^2}$, ${i=0,1,\ldots ,n}$. For any parameter ${x\in\hat{\Om}}$ and B-spline bases given by~\eqref{eq.BasisFunBip}, it reads
\begin{align}
	C(x)=\sum_{i=0}^{n} \B{i}{p}(x)\.P_i \,.
	\label{eq.nonrational}
\end{align}
Let us consider quadratic bases (i.e., ${p=2}$) spanned over the knot vector ${\Xi=\{0,0,0,\frac{1}{5},\frac{2}{5},\frac{3}{5},\frac{4}{5},1,1,1\}}$.
\autoref{fig:Tschirn} shows the midline of the Tschirnhausen beam and its corresponding control points listed in \autoref{tab:TschirnCP}. 
We use the algebraic formula of the Tschirnhausen curve (see, e.g.,~\cite{ Farouki2008, Hosseini2018}) and, then, generate the control points by a curve fitting process~\cite{TheNURBSBook}.

\begin{table}[!h]\centering
	\caption{Control points of the Tschirnhausen beam.}
	\label{tab:TschirnCP}
    \small
	\begin{tabular}{@{}cccccccc@{}}
		\toprule
		Control points & $P_0$  & $P_1$  & $P_2$  & $P_3$  & $P_4$  & $P_5$  & $P_6$  \\ \midrule
		$x$-component  & 0.0000 & 0.0122 & 0.2121 & 0.4196 & 0.6208 & 0.8099 & 0.9000 \\
		$y$-component  & 0.0000 & 0.1137 & 0.2065 & 0.1945 & 0.1385 & 0.0519 & 0.0000 \\ \bottomrule
	\end{tabular}
\end{table}

\subsubsection{Eigenfrequency analysis}

We use the Timoshenko curved beams formulation; but details are omitted for the sake of brevity (see, e.g.,~\cite{Luu2014,Hosseini20182} for complete formulation).
In this context, each nodal component of the beam has three degrees of freedom, namely, tangential and normal displacements of the midline, $\ur$ and $\wr$, respectively, and transverse rotation of the cross section $\vphi$.
We assume our computational domain $\Om$ as a mapping from the rectilinear parameter domain $\hat{\Om}$ onto the curved geometry.
Let us introduce ${\ui:=\{\ur,\wr,\vphi\}}$ with ${\ui:\Om\rightarrow\R^3}$ as the vector of degrees of freedom, and ${\vi\in\big(H^1_0(\Om)\big)^3}$ as our test space.
We write the eigenproblem associated with the frequency analysis of curved beams in the weak form as follows:
Find ${\lm\in\R^+}$ and ${\ui\in\big(H^1_0(\Om)\big)^3}$, such that for all ${\vi\in\big(H^1_0(\Om)\big)^3}$,
\begin{align}
\Ar(\vi,\ui) = \lm\.\Br(\vi,\ui)\,,
\label{eq:TschirnWeak}
\end{align}
where $\lm^{1/2}$ is the eigenfrequency. We define the bilinear forms
\begin{align}
\Ar(\vi,\ui) &:= \int_\Om \eps(\vi)\cdot\sg(\ui) \,d\Om\,, \label{eq:BilinearA}\\
\Br(\vi,\ui) &:= \int_\Om \vi\cdot \Gm\ui \,d\Om\,, \label{eq:BilinearB}
\end{align}
with
\begin{align}
\eps(\ui) &= \big\{ \ur'-\kappa\wr \,,\, \wr'+\kappa\ur-\vphi \,,\, \vphi' \big\} \,, \label{eq:StrainField}\\
\sg(\ui) &= {\rm diag}\.\big\{ EA \,,\, k_sGA \,,\, EI \big\} \. \eps(\ui) \,,\\
\Gm &= {\rm diag}\.\big\{ \rho A \,,\, \rho A \,,\, \rho I \big\} \,,
\end{align}
where $\kappa$ is the curvature of the midline, $\rho$ is the material density, $E$ and $G$ are the elastic and shear modulus of the material of the beam, $A$ is cross-section area, $I$ is the second moment of inertia of cross-section, and $k_s$ is the shear correction factor.

In the IGA framework, it is a common practice to use the same basis functions that define the geometry to discretize the solution fields $\ur$, $\wr$, and $\vphi$.
Let us introduce ${\hat{\Vi}^h:=\big(\Sr^p_c(U)\big)^3}$ as our discrete spaces in the parametric domain.
We define the discrete space in the physical domain
\begin{align}
    \Vi^h := \Big\{\vi^h\in\big(H^1_0(\Om)\big)^3:\iota(\vi^h)\in\hat{\Vi}^h\Big\}\,,
\end{align}
where ${\iota(\cdot)}$ is the pullback mapping to the parameter space (see, e.g.,~\cite{Garcia2019}).
Thus, we write the discrete form of the eigenproblem~\eqref{eq:TschirnWeak} as follows:
Find ${\lm^h\in\R^+}$ and ${\ui^h\in \Vi^h}$, such that for all ${\vi^h\in \Vi^h}$,
\begin{align}
\Ar(\vi^h,\ui^h) = \lm^h\.\Br(\vi^h,\ui^h)\,,
\label{eq:TschirnWeakDiscrete}
\end{align}
that results in the same matrix form as~\eqref{eq:EigProblem}, while the total number of degrees of freedom $\N$ is increased to~${3\.(n-2)}$.
When constructing the stiffness and mass matrices given by the bilinear forms~\eqref{eq:BilinearA} and~\eqref{eq:BilinearB}, we integrate over the curved domain $\Om$ noting that the derivatives in~\eqref{eq:StrainField} are also expressed in the curvilinear coordinate system.
Let 
\begin{align}
J:=\norm{\dfrac{d\.C(x)}{dx}}
\end{align}
be the \textit{Jacobian}, which is the Euclidean norm of the derivative of the geometry, we write ${d\Om=J\.d\hat{\Om}}$ and ${(\cdot)'=\frac{1}{J}\frac{d(\cdot)}{dx}}$ to compute the derivatives and integrals in the rectilinear parameter domain $\hat{\Om}$. In this notation, we compute the curvature of the geometry as follows:
\begin{align}
\kappa:=\dfrac{1}{J}\norm{\dfrac{d\.C(x)}{dx}\times\dfrac{d^2C(x)}{dx^2}}.
\end{align}

\begin{figure}[!h]\centering
	\includegraphics{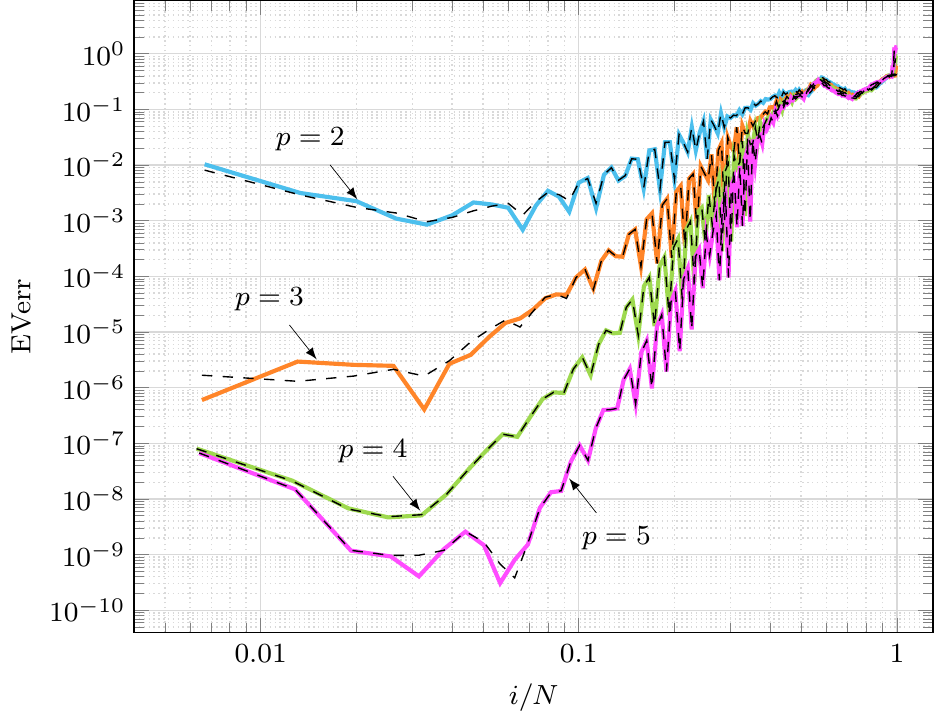}
	\caption{Eigenfrequency error analysis of the Tschirnhausen beam discretized by ${\ne = 50}$ uniform elements and different polynomial degrees. The difference between results of the optimal and EWG rules (solid and dashed lines, respectively) is because of the variation of the curvature.}
	\label{fig:EVErrTschirn}	
\end{figure}

In the following, we compare the eigenvalue errors of the Tschirnhausen curved beam when constructing the system matrices using the optimal and EWG quadrature rules (see \autoref{fig:EVErrTschirn}).
We consider ${\ne = 50}$ uniform elements and polynomial degrees ${p=2,3,4,5}$ to discretize the respective eigenproblem. 
The difference between results of the optimal and EWG rules is because of the variation of the curvature and the inclusion of the geometric terms in integrals~\eqref{eq:BilinearA} and~\eqref{eq:BilinearB}.
Nevertheless, neither of these rules are exact in this case because they are constructed based on~\eqref{eq:Stiffness} and~\eqref{eq:Mass} without taking into account the geometric terms in the integration (c.f.~\cite[Section~3.2]{hughes_efficient_2010} and~\cite[Section~2]{ auricchio_simple_2012}).
Finally, we note that the low-frequency region of the spectrum of this curved beam shows relatively high frequency errors under the employment of typical IGA discretizations (see, e.g.,~\cite{Hosseini2018,Hosseini2020}).

\section{Conclusions}
\label{sec:Conclusions}

We have developed and validated a novel machine-learning method for finding optimal quadrature rules on B-spline spaces. This method has shown general convergence to the optimal rule on uniform elements with arbitrary degrees and continuities, with a good computational scalability. We have also employed the method to construct system matrices for practical IGA problems to assess its accuracy. In particular, we solved the eigenproblem of the Laplace operator and eigenfrequency problem of a planar curved beam. We showed that the results are equivalent to element-wise Gaussian integration, but requiring up to 44\% fewer integration points. These savings are magnified when performing integrals over high-dimensional basis functions, reaching 68\% in 2D spaces, and up to 82\% in 3D spaces. Finally, even if we have demonstrated that the method works for some partitions with non-uniform elements, as a future line of research we will pursue a strategy that is generally applicable in this scenario.


\section*{Acknowledgements}
This work has been funded by
the Euskampus Foundation through the ORLEG-IA project in the Misiones Euskampus 2.0 program,
the ``BCAM Severo Ochoa'' accreditation of excellence CEX2021-001142-S / MICIN / AEI / 10.13039/501100011033, 
the projects of the Spanish Ministry of Science and Innovation
PID2019-108111RB-I00,
PDC2021-121093-I00,
TED2021-132783B-I00,
and the Basque Government through the BERC 2022-2025 program and the Consolidated Research Group MATHMODE (IT1456-22) given by the Department of Education. 
T.~Teijeiro is supported by the grant RYC2021-032853-I funded by MCIN/AEI/ 10.13039/501100011033 and by the EU NextGenerationEU/PRTR.


\bibliography{bibliography}

\begin{thebibliography}{39}
\expandafter\ifx\csname natexlab\endcsname\relax\def\natexlab#1{#1}\fi
\providecommand{\url}[1]{\texttt{#1}}
\providecommand{\href}[2]{#2}
\providecommand{\path}[1]{#1}
\providecommand{\DOIprefix}{doi:}
\providecommand{\ArXivprefix}{arXiv:}
\providecommand{\URLprefix}{URL: }
\providecommand{\Pubmedprefix}{pmid:}
\providecommand{\doi}[1]{\href{http://dx.doi.org/#1}{\path{#1}}}
\providecommand{\Pubmed}[1]{\href{pmid:#1}{\path{#1}}}
\providecommand{\bibinfo}[2]{#2}
\ifx\xfnm\relax \def\xfnm[#1]{\unskip,\space#1}\fi
\bibitem[{Hughes et~al.(2005)Hughes, Cottrell, and Bazilevs}]{Hughes2005}
\bibinfo{author}{T.~J.~R. Hughes}, \bibinfo{author}{J.~A. Cottrell},
  \bibinfo{author}{Y.~Bazilevs},
\newblock \bibinfo{title}{{Isogeometric analysis: CAD, finite elements, NURBS,
  exact geometry and mesh refinement}},
\newblock \bibinfo{journal}{Computer Methods in Applied Mechanics and
  Engineering} \bibinfo{volume}{194} (\bibinfo{year}{2005})
  \bibinfo{pages}{4135--4195}. \DOIprefix\doi{10.1016/j.cma.2004.10.008}.
\bibitem[{Calabr{\`{o}} et~al.(2017)Calabr{\`{o}}, Sangalli, and
  Tani}]{Calabro2017}
\bibinfo{author}{F.~Calabr{\`{o}}}, \bibinfo{author}{G.~Sangalli},
  \bibinfo{author}{M.~Tani},
\newblock \bibinfo{title}{Fast formation of isogeometric {Galerkin} matrices by
  weighted quadrature},
\newblock \bibinfo{journal}{Computer Methods in Applied Mechanics and
  Engineering} \bibinfo{volume}{316} (\bibinfo{year}{2017})
  \bibinfo{pages}{606--622}. \DOIprefix\doi{10.1016/j.cma.2016.09.013}.
\bibitem[{Barto{\v{n}} et~al.(2020)Barto{\v{n}}, Puzyrev, Deng, and
  Calo}]{Barton2020}
\bibinfo{author}{M.~Barto{\v{n}}}, \bibinfo{author}{V.~Puzyrev},
  \bibinfo{author}{Q.~Deng}, \bibinfo{author}{V.~Calo},
\newblock \bibinfo{title}{Efficient mass and stiffness matrix assembly via
  weighted {Gaussian} quadrature rules for {B-splines}},
\newblock \bibinfo{journal}{Journal of Computational and Applied Mathematics}
  \bibinfo{volume}{371} (\bibinfo{year}{2020}) \bibinfo{pages}{112626}.
  \DOIprefix\doi{10.1016/j.cam.2019.112626}.
\bibitem[{Barto{\v n} and Calo(2016{\natexlab{a}})}]{barton_gaussian_2016}
\bibinfo{author}{M.~Barto{\v n}}, \bibinfo{author}{V.~M. Calo},
\newblock \bibinfo{title}{Gaussian quadrature for splines via homotopy
  continuation: Rules for ${C^2}$ cubic splines},
\newblock \bibinfo{journal}{Journal of Computational and Applied Mathematics}
  \bibinfo{volume}{296} (\bibinfo{year}{2016}{\natexlab{a}})
  \bibinfo{pages}{709--723}. \DOIprefix\doi{10.1016/j.cam.2015.09.036}.
\bibitem[{Barto{\v n} and
  Calo(2016{\natexlab{b}})}]{bartonOptimalQuadratureRules2016}
\bibinfo{author}{M.~Barto{\v n}}, \bibinfo{author}{V.~M. Calo},
\newblock \bibinfo{title}{Optimal quadrature rules for odd-degree spline spaces
  and their application to tensor-product-based isogeometric analysis},
\newblock \bibinfo{journal}{Computer Methods in Applied Mechanics and
  Engineering} \bibinfo{volume}{305} (\bibinfo{year}{2016}{\natexlab{b}})
  \bibinfo{pages}{217--240}. \DOIprefix\doi{10.1016/j.cma.2016.02.034}.
\bibitem[{Barto{\v n} and Calo(2017)}]{Barton2017}
\bibinfo{author}{M.~Barto{\v n}}, \bibinfo{author}{V.~M. Calo},
\newblock \bibinfo{title}{{Gauss--Galerkin} quadrature rules for quadratic and
  cubic spline spaces and their application to isogeometric analysis},
\newblock \bibinfo{journal}{Computer-Aided Design} \bibinfo{volume}{82}
  (\bibinfo{year}{2017}) \bibinfo{pages}{57--67}.
  \DOIprefix\doi{10.1016/j.cad.2016.07.003}.
\bibitem[{Sommese and Wampler(2005)}]{Sommese2005}
\bibinfo{author}{A.~J. Sommese}, \bibinfo{author}{C.~W. Wampler},
  \bibinfo{title}{The Numerical Solution of Systems of Polynomials Arising in
  Engineering and Science}, \bibinfo{publisher}{World Scientific},
  \bibinfo{address}{Singapore}, \bibinfo{year}{2005}.
\bibitem[{Hughes et~al.(2010)Hughes, Reali, and
  Sangalli}]{hughes_efficient_2010}
\bibinfo{author}{T.~J.~R. Hughes}, \bibinfo{author}{A.~Reali},
  \bibinfo{author}{G.~Sangalli},
\newblock \bibinfo{title}{Efficient quadrature for {NURBS}-based isogeometric
  analysis},
\newblock \bibinfo{journal}{Computer Methods in Applied Mechanics and
  Engineering} \bibinfo{volume}{199} (\bibinfo{year}{2010})
  \bibinfo{pages}{301--313}. \DOIprefix\doi{10.1016/j.cma.2008.12.004}.
\bibitem[{Auricchio et~al.(2012)Auricchio, Calabr{\`{o}}, Hughes, Reali, and
  Sangalli}]{auricchio_simple_2012}
\bibinfo{author}{F.~Auricchio}, \bibinfo{author}{F.~Calabr{\`{o}}},
  \bibinfo{author}{T.~J.~R. Hughes}, \bibinfo{author}{A.~Reali},
  \bibinfo{author}{G.~Sangalli},
\newblock \bibinfo{title}{A simple algorithm for obtaining nearly optimal
  quadrature rules for {NURBS}-based isogeometric analysis},
\newblock \bibinfo{journal}{Computer Methods in Applied Mechanics and
  Engineering} \bibinfo{volume}{249-252} (\bibinfo{year}{2012})
  \bibinfo{pages}{15--27}. \DOIprefix\doi{10.1016/j.cma.2012.04.014}.
\bibitem[{Schillinger et~al.(2014)Schillinger, Hossain, and
  Hughes}]{Schillinger2014}
\bibinfo{author}{D.~Schillinger}, \bibinfo{author}{S.~J. Hossain},
  \bibinfo{author}{T.~J.~R. Hughes},
\newblock \bibinfo{title}{Reduced b{\'{e}}zier element quadrature rules for
  quadratic and cubic splines in isogeometric analysis},
\newblock \bibinfo{journal}{Computer Methods in Applied Mechanics and
  Engineering} \bibinfo{volume}{277} (\bibinfo{year}{2014})
  \bibinfo{pages}{1--45}. \DOIprefix\doi{10.1016/j.cma.2014.04.008}.
\bibitem[{Hiemstra et~al.(2017)Hiemstra, Calabr{\`{o}}, Schillinger, and
  Hughes}]{Hiemstra2017}
\bibinfo{author}{R.~R. Hiemstra}, \bibinfo{author}{F.~Calabr{\`{o}}},
  \bibinfo{author}{D.~Schillinger}, \bibinfo{author}{T.~J.~R. Hughes},
\newblock \bibinfo{title}{Optimal and reduced quadrature rules for tensor
  product and hierarchically refined splines in isogeometric analysis},
\newblock \bibinfo{journal}{Computer Methods in Applied Mechanics and
  Engineering} \bibinfo{volume}{316} (\bibinfo{year}{2017})
  \bibinfo{pages}{966--1004}. \DOIprefix\doi{10.1016/j.cma.2016.10.049}.
\bibitem[{Barendrecht et~al.(2018)Barendrecht, Barto{\v{n}}, and
  Kosinka}]{Barendrecht2018}
\bibinfo{author}{P.~J. Barendrecht}, \bibinfo{author}{M.~Barto{\v{n}}},
  \bibinfo{author}{J.~Kosinka},
\newblock \bibinfo{title}{Efficient quadrature rules for subdivision surfaces
  in isogeometric analysis},
\newblock \bibinfo{journal}{Computer Methods in Applied Mechanics and
  Engineering} \bibinfo{volume}{340} (\bibinfo{year}{2018})
  \bibinfo{pages}{1--23}. \DOIprefix\doi{10.1016/j.cma.2018.05.017}.
\bibitem[{Zou et~al.(2021)Zou, Hughes, Scott, Sauer, and Savitha}]{Zou2021}
\bibinfo{author}{Z.~Zou}, \bibinfo{author}{T.~J.~R. Hughes},
  \bibinfo{author}{M.~A. Scott}, \bibinfo{author}{R.~A. Sauer},
  \bibinfo{author}{E.~J. Savitha},
\newblock \bibinfo{title}{Galerkin formulations of isogeometric shell analysis:
  Alleviating locking with {Greville} quadratures and higher-order elements},
\newblock \bibinfo{journal}{Computer Methods in Applied Mechanics and
  Engineering} \bibinfo{volume}{380} (\bibinfo{year}{2021})
  \bibinfo{pages}{113757}. \DOIprefix\doi{10.1016/j.cma.2021.113757}.
\bibitem[{Giannelli et~al.(2022)Giannelli, Kandu{\v{c}}, Martinelli, Sangalli,
  and Tani}]{Giannelli2022}
\bibinfo{author}{C.~Giannelli}, \bibinfo{author}{T.~Kandu{\v{c}}},
  \bibinfo{author}{M.~Martinelli}, \bibinfo{author}{G.~Sangalli},
  \bibinfo{author}{M.~Tani},
\newblock \bibinfo{title}{Weighted quadrature for hierarchical {B-splines}},
\newblock \bibinfo{journal}{Computer Methods in Applied Mechanics and
  Engineering} \bibinfo{volume}{400} (\bibinfo{year}{2022})
  \bibinfo{pages}{115465}. \DOIprefix\doi{10.1016/j.cma.2022.115465}.
\bibitem[{Mon(2004)}]{MonteCarloIntegration2004}
\bibinfo{title}{Monte {Carlo} integration {{I}}},
\newblock in: \bibinfo{booktitle}{Physically {{Based Rendering}}},
  \bibinfo{publisher}{{Elsevier}}, \bibinfo{year}{2004}, pp.
  \bibinfo{pages}{631--660}. \DOIprefix\doi{10.1016/B978-012553180-1/50016-8}.
\bibitem[{Kanagawa et~al.(2016)Kanagawa, Sriperumbudur, and
  Fukumizu}]{NIPS2016_81c650ca}
\bibinfo{author}{M.~Kanagawa}, \bibinfo{author}{B.~K. Sriperumbudur},
  \bibinfo{author}{K.~Fukumizu},
\newblock \bibinfo{title}{Convergence guarantees for kernel-based quadrature
  rules in misspecified settings},
\newblock in: \bibinfo{editor}{D.~Lee}, \bibinfo{editor}{M.~Sugiyama},
  \bibinfo{editor}{U.~Luxburg}, \bibinfo{editor}{I.~Guyon},
  \bibinfo{editor}{R.~Garnett} (Eds.), \bibinfo{booktitle}{Advances in Neural
  Information Processing Systems}, volume~\bibinfo{volume}{29},
  \bibinfo{publisher}{Curran Associates, Inc.}, \bibinfo{year}{2016}.
  \DOIprefix\doi{10.5555/3157382.3157466}.
\bibitem[{Munkhoeva et~al.(2018)Munkhoeva, Kapushev, Burnaev, and
  Oseledets}]{NEURIPS2018_6e923226}
\bibinfo{author}{M.~Munkhoeva}, \bibinfo{author}{Y.~Kapushev},
  \bibinfo{author}{E.~Burnaev}, \bibinfo{author}{I.~Oseledets},
\newblock \bibinfo{title}{Quadrature-based features for kernel approximation},
\newblock in: \bibinfo{editor}{S.~Bengio}, \bibinfo{editor}{H.~Wallach},
  \bibinfo{editor}{H.~Larochelle}, \bibinfo{editor}{K.~Grauman},
  \bibinfo{editor}{N.~Cesa-Bianchi}, \bibinfo{editor}{R.~Garnett} (Eds.),
  \bibinfo{booktitle}{Advances in Neural Information Processing Systems},
  volume~\bibinfo{volume}{31}, \bibinfo{publisher}{Curran Associates, Inc.},
  \bibinfo{year}{2018}. \DOIprefix\doi{10.5555/3327546.3327589}.
\bibitem[{O'Hagan(1991)}]{ohaganBayesHermiteQuadrature1991}
\bibinfo{author}{A.~O'Hagan},
\newblock \bibinfo{title}{Bayes\textendash{{Hermite}} quadrature},
\newblock \bibinfo{journal}{Journal of Statistical Planning and Inference}
  \bibinfo{volume}{29} (\bibinfo{year}{1991}) \bibinfo{pages}{245--260}.
  \DOIprefix\doi{10.1016/0378-3758(91)90002-V}.
\bibitem[{Karvonen and Sarkka(2017)}]{karvonenClassicalQuadratureRules2017}
\bibinfo{author}{T.~Karvonen}, \bibinfo{author}{S.~Sarkka},
\newblock \bibinfo{title}{Classical quadrature rules via {{Gaussian}}
  processes},
\newblock in: \bibinfo{booktitle}{2017 {{IEEE}} 27th {{International Workshop}}
  on {{Machine Learning}} for {{Signal Processing}} ({{MLSP}})},
  \bibinfo{publisher}{{IEEE}}, \bibinfo{address}{{Tokyo}},
  \bibinfo{year}{2017}, pp. \bibinfo{pages}{1--6}.
  \DOIprefix\doi{10.1109/MLSP.2017.8168195}.
\bibitem[{Kanagawa and Hennig(2019)}]{NEURIPS2019_165a59f7}
\bibinfo{author}{M.~Kanagawa}, \bibinfo{author}{P.~Hennig},
\newblock \bibinfo{title}{Convergence guarantees for adaptive bayesian
  quadrature methods},
\newblock in: \bibinfo{editor}{H.~Wallach}, \bibinfo{editor}{H.~Larochelle},
  \bibinfo{editor}{A.~Beygelzimer}, \bibinfo{editor}{F.~d\textquotesingle
  Alch\'{e}-Buc}, \bibinfo{editor}{E.~Fox}, \bibinfo{editor}{R.~Garnett}
  (Eds.), \bibinfo{booktitle}{Advances in Neural Information Processing
  Systems}, volume~\bibinfo{volume}{32}, \bibinfo{publisher}{Curran Associates,
  Inc.}, \bibinfo{year}{2019}. \DOIprefix\doi{10.5555/3454287.3454847}.
\bibitem[{Rivera et~al.(2022)Rivera, Taylor, Omella, and
  Pardo}]{riveraQuadratureRulesSolving2022}
\bibinfo{author}{J.~A. Rivera}, \bibinfo{author}{J.~M. Taylor},
  \bibinfo{author}{{\'A}.~J. Omella}, \bibinfo{author}{D.~Pardo},
\newblock \bibinfo{title}{On quadrature rules for solving {{Partial
  Differential Equations}} using {{Neural Networks}}},
\newblock \bibinfo{journal}{Computer Methods in Applied Mechanics and
  Engineering} \bibinfo{volume}{393} (\bibinfo{year}{2022})
  \bibinfo{pages}{114710}. \DOIprefix\doi{10.1016/j.cma.2022.114710}.
\bibitem[{Narkhede et~al.(2022)Narkhede, Bartakke, and
  Sutaone}]{narkhedeReviewWeightInitialization2022}
\bibinfo{author}{M.~V. Narkhede}, \bibinfo{author}{P.~P. Bartakke},
  \bibinfo{author}{M.~S. Sutaone},
\newblock \bibinfo{title}{A review on weight initialization strategies for
  neural networks},
\newblock \bibinfo{journal}{Artificial Intelligence Review}
  \bibinfo{volume}{55} (\bibinfo{year}{2022}) \bibinfo{pages}{291--322}.
  \DOIprefix\doi{10.1007/s10462-021-10033-z}.
\bibitem[{Bellman(1966)}]{bellmanDynamicProgramming1966}
\bibinfo{author}{R.~Bellman},
\newblock \bibinfo{title}{Dynamic {{Programming}}},
\newblock \bibinfo{journal}{Science} \bibinfo{volume}{153}
  (\bibinfo{year}{1966}) \bibinfo{pages}{34--37}.
  \DOIprefix\doi{10.1126/science.153.3731.34}.
\bibitem[{Piegl and Tiller(1996)}]{TheNURBSBook}
\bibinfo{author}{L.~Piegl}, \bibinfo{author}{W.~Tiller}, \bibinfo{title}{The
  {NURBS} book}, \bibinfo{publisher}{Springer-Verlag}, \bibinfo{year}{1996}.
\bibitem[{Hashemian et~al.(2022)Hashemian, Garcia, Pardo, and
  Calo}]{Hashemian2022}
\bibinfo{author}{A.~Hashemian}, \bibinfo{author}{D.~Garcia},
  \bibinfo{author}{D.~Pardo}, \bibinfo{author}{V.~M. Calo},
\newblock \bibinfo{title}{Refined isogeometric analysis of quadratic eigenvalue
  problems},
\newblock \bibinfo{journal}{Computer Methods in Applied Mechanics and
  Engineering} \bibinfo{volume}{399} (\bibinfo{year}{2022})
  \bibinfo{pages}{115327}. \DOIprefix\doi{10.1016/j.cma.2022.115327}.
\bibitem[{Zaheer et~al.(2018)Zaheer, Reddi, Sachan, Kale, and
  Kumar}]{NEURIPS18_adaptive}
\bibinfo{author}{M.~Zaheer}, \bibinfo{author}{S.~Reddi},
  \bibinfo{author}{D.~Sachan}, \bibinfo{author}{S.~Kale},
  \bibinfo{author}{S.~Kumar},
\newblock \bibinfo{title}{Adaptive methods for nonconvex optimization},
\newblock in: \bibinfo{editor}{S.~Bengio}, \bibinfo{editor}{H.~Wallach},
  \bibinfo{editor}{H.~Larochelle}, \bibinfo{editor}{K.~Grauman},
  \bibinfo{editor}{N.~Cesa-Bianchi}, \bibinfo{editor}{R.~Garnett} (Eds.),
  \bibinfo{booktitle}{Advances in Neural Information Processing Systems},
  volume~\bibinfo{volume}{31}, \bibinfo{publisher}{Curran Associates, Inc.},
  \bibinfo{year}{2018}.
\bibitem[{Bradbury et~al.(2018)Bradbury, Frostig, Hawkins, Johnson, Leary,
  Maclaurin, Necula, Paszke, Vander{P}las, Wanderman-{M}ilne, and
  Zhang}]{jax2018github}
\bibinfo{author}{J.~Bradbury}, \bibinfo{author}{R.~Frostig},
  \bibinfo{author}{P.~Hawkins}, \bibinfo{author}{M.~J. Johnson},
  \bibinfo{author}{C.~Leary}, \bibinfo{author}{D.~Maclaurin},
  \bibinfo{author}{G.~Necula}, \bibinfo{author}{A.~Paszke},
  \bibinfo{author}{J.~Vander{P}las}, \bibinfo{author}{S.~Wanderman-{M}ilne},
  \bibinfo{author}{Q.~Zhang}, \bibinfo{title}{{JAX}: composable transformations
  of {P}ython+{N}um{P}y programs}, \bibinfo{year}{2018}. \URLprefix
  \url{http://github.com/google/jax}.
\bibitem[{Babuschkin et~al.(2020)Babuschkin, Baumli, Bell, Bhupatiraju, Bruce,
  Buchlovsky, Budden, Cai, Clark, Danihelka, Dedieu, Fantacci, Godwin, Jones,
  Hemsley, Hennigan, Hessel, Hou, Kapturowski, Keck, Kemaev, King, Kunesch,
  Martens, Merzic, Mikulik, Norman, Papamakarios, Quan, Ring, Ruiz, Sanchez,
  Schneider, Sezener, Spencer, Srinivasan, Stokowiec, Wang, Zhou, and
  Viola}]{deepmind2020jax}
\bibinfo{author}{I.~Babuschkin}, \bibinfo{author}{K.~Baumli},
  \bibinfo{author}{A.~Bell}, \bibinfo{author}{S.~Bhupatiraju},
  \bibinfo{author}{J.~Bruce}, \bibinfo{author}{P.~Buchlovsky},
  \bibinfo{author}{D.~Budden}, \bibinfo{author}{T.~Cai},
  \bibinfo{author}{A.~Clark}, \bibinfo{author}{I.~Danihelka},
  \bibinfo{author}{A.~Dedieu}, \bibinfo{author}{C.~Fantacci},
  \bibinfo{author}{J.~Godwin}, \bibinfo{author}{C.~Jones},
  \bibinfo{author}{R.~Hemsley}, \bibinfo{author}{T.~Hennigan},
  \bibinfo{author}{M.~Hessel}, \bibinfo{author}{S.~Hou},
  \bibinfo{author}{S.~Kapturowski}, \bibinfo{author}{T.~Keck},
  \bibinfo{author}{I.~Kemaev}, \bibinfo{author}{M.~King},
  \bibinfo{author}{M.~Kunesch}, \bibinfo{author}{L.~Martens},
  \bibinfo{author}{H.~Merzic}, \bibinfo{author}{V.~Mikulik},
  \bibinfo{author}{T.~Norman}, \bibinfo{author}{G.~Papamakarios},
  \bibinfo{author}{J.~Quan}, \bibinfo{author}{R.~Ring},
  \bibinfo{author}{F.~Ruiz}, \bibinfo{author}{A.~Sanchez},
  \bibinfo{author}{R.~Schneider}, \bibinfo{author}{E.~Sezener},
  \bibinfo{author}{S.~Spencer}, \bibinfo{author}{S.~Srinivasan},
  \bibinfo{author}{W.~Stokowiec}, \bibinfo{author}{L.~Wang},
  \bibinfo{author}{G.~Zhou}, \bibinfo{author}{F.~Viola}, \bibinfo{title}{The
  {D}eep{M}ind {JAX} {E}cosystem}, \bibinfo{year}{2020}. \URLprefix
  \url{http://github.com/deepmind}.
\bibitem[{Hashemian et~al.(2021)Hashemian, Pardo, and Calo}]{Hashemian2021}
\bibinfo{author}{A.~Hashemian}, \bibinfo{author}{D.~Pardo},
  \bibinfo{author}{V.~M. Calo},
\newblock \bibinfo{title}{Refined isogeometric analysis for generalized
  {Hermitian} eigenproblems},
\newblock \bibinfo{journal}{Computer Methods in Applied Mechanics and
  Engineering} \bibinfo{volume}{381} (\bibinfo{year}{2021})
  \bibinfo{pages}{113823}. \DOIprefix\doi{10.1016/j.cma.2021.113823}.
\bibitem[{Puzyrev et~al.(2017)Puzyrev, Deng, and Calo}]{Puzyrev2017}
\bibinfo{author}{V.~Puzyrev}, \bibinfo{author}{Q.~Deng},
  \bibinfo{author}{V.~Calo},
\newblock \bibinfo{title}{Dispersion-optimized quadrature rules for
  isogeometric analysis: Modified inner products, their dispersion properties,
  and optimally blended schemes},
\newblock \bibinfo{journal}{Computer Methods in Applied Mechanics and
  Engineering} \bibinfo{volume}{320} (\bibinfo{year}{2017})
  \bibinfo{pages}{421--443}. \DOIprefix\doi{10.1016/j.cma.2017.03.029}.
\bibitem[{Garcia et~al.(2019)Garcia, Pardo, and Calo}]{Garcia2019}
\bibinfo{author}{D.~Garcia}, \bibinfo{author}{D.~Pardo}, \bibinfo{author}{V.~M.
  Calo},
\newblock \bibinfo{title}{Refined isogeometric analysis for fluid mechanics and
  electromagnetics},
\newblock \bibinfo{journal}{Computer Methods in Applied Mechanics and
  Engineering} \bibinfo{volume}{356} (\bibinfo{year}{2019})
  \bibinfo{pages}{598--628}. \DOIprefix\doi{10.1016/j.cma.2019.06.011}.
\bibitem[{Hashemian et~al.(2021)Hashemian, Garcia, Rivera, and
  Pardo}]{Hashemian20212}
\bibinfo{author}{A.~Hashemian}, \bibinfo{author}{D.~Garcia},
  \bibinfo{author}{J.~A. Rivera}, \bibinfo{author}{D.~Pardo},
\newblock \bibinfo{title}{Massive database generation for {2.5D} borehole
  electromagnetic measurements using refined isogeometric analysis},
\newblock \bibinfo{journal}{Computers {\&} Geosciences} \bibinfo{volume}{155}
  (\bibinfo{year}{2021}) \bibinfo{pages}{104808}.
  \DOIprefix\doi{10.1016/j.cageo.2021.104808}.
\bibitem[{Garcia et~al.(2017)Garcia, Pardo, Dalcin, Paszy\'nski, Collier, and
  Calo}]{Garcia2017}
\bibinfo{author}{D.~Garcia}, \bibinfo{author}{D.~Pardo},
  \bibinfo{author}{L.~Dalcin}, \bibinfo{author}{M.~Paszy\'nski},
  \bibinfo{author}{N.~Collier}, \bibinfo{author}{V.~M. Calo},
\newblock \bibinfo{title}{The value of continuity: Refined isogeometric
  analysis and fast direct solvers},
\newblock \bibinfo{journal}{Computer Methods in Applied Mechanics and
  Engineering} \bibinfo{volume}{316} (\bibinfo{year}{2017})
  \bibinfo{pages}{586--605}. \DOIprefix\doi{10.1016/j.cma.2016.08.017}.
\bibitem[{Puzyrev et~al.(2018)Puzyrev, Deng, and Calo}]{Puzyrev2018}
\bibinfo{author}{V.~Puzyrev}, \bibinfo{author}{Q.~Deng},
  \bibinfo{author}{V.~Calo},
\newblock \bibinfo{title}{Spectral approximation properties of isogeometric
  analysis with variable continuity},
\newblock \bibinfo{journal}{Computer Methods in Applied Mechanics and
  Engineering} \bibinfo{volume}{334} (\bibinfo{year}{2018})
  \bibinfo{pages}{22--39}. \DOIprefix\doi{10.1016/j.cma.2018.01.042}.
\bibitem[{Farouki(2008)}]{Farouki2008}
\bibinfo{author}{R.~T. Farouki}, \bibinfo{title}{Pythagorean-Hodograph Curves:
  Algebra and Geometry Inseparable}, \bibinfo{publisher}{Springer Berlin
  Heidelberg}, \bibinfo{year}{2008}.
\bibitem[{Hosseini et~al.(2018)Hosseini, Hashemian, and Reali}]{Hosseini2018}
\bibinfo{author}{S.~F. Hosseini}, \bibinfo{author}{A.~Hashemian},
  \bibinfo{author}{A.~Reali},
\newblock \bibinfo{title}{On the application of curve reparameterization in
  isogeometric vibration analysis of free-from curved beams},
\newblock \bibinfo{journal}{Computers {\&} Structures} \bibinfo{volume}{209}
  (\bibinfo{year}{2018}) \bibinfo{pages}{117--129}.
  \DOIprefix\doi{10.1016/j.compstruc.2018.08.009}.
\bibitem[{Luu et~al.(2014)Luu, Kim, and Lee}]{Luu2014}
\bibinfo{author}{A.-T. Luu}, \bibinfo{author}{N.-I. Kim},
  \bibinfo{author}{J.~Lee},
\newblock \bibinfo{title}{Isogeometric vibration analysis of free-form
  {Timoshenko} curved beams},
\newblock \bibinfo{journal}{Meccanica} \bibinfo{volume}{50}
  (\bibinfo{year}{2014}) \bibinfo{pages}{169--187}.
  \DOIprefix\doi{10.1007/s11012-014-0062-3}.
\bibitem[{Hosseini et~al.(2018)Hosseini, Hashemian, Moetakef-Imani, and
  Hadidimoud}]{Hosseini20182}
\bibinfo{author}{S.~F. Hosseini}, \bibinfo{author}{A.~Hashemian},
  \bibinfo{author}{B.~Moetakef-Imani}, \bibinfo{author}{S.~Hadidimoud},
\newblock \bibinfo{title}{Isogeometric analysis of free-form {Timoshenko}
  curved beams including the nonlinear effects of large deformations},
\newblock \bibinfo{journal}{Acta Mechanica Sinica} \bibinfo{volume}{34}
  (\bibinfo{year}{2018}) \bibinfo{pages}{728--743}.
  \DOIprefix\doi{10.1007/s10409-018-0753-4}.
\bibitem[{Hosseini et~al.(2020)Hosseini, Hashemian, and Reali}]{Hosseini2020}
\bibinfo{author}{S.~F. Hosseini}, \bibinfo{author}{A.~Hashemian},
  \bibinfo{author}{A.~Reali},
\newblock \bibinfo{title}{Studies on knot placement techniques for the geometry
  construction and the accurate simulation of isogeometric spatial curved
  beams},
\newblock \bibinfo{journal}{Computer Methods in Applied Mechanics and
  Engineering} \bibinfo{volume}{360} (\bibinfo{year}{2020})
  \bibinfo{pages}{112705}. \DOIprefix\doi{10.1016/j.cma.2019.112705}.

\end{thebibliography}

\end{document}